\documentclass[preprint,review,11pt]{elsarticle}

\usepackage{amssymb}
\usepackage{amsmath}
\usepackage{amsthm}
\usepackage{esint}
\newtheorem{theorem}{Theorem}[section]

\numberwithin{equation}{section}

\newtheorem{corollary}[theorem]{Corollary}

\newtheorem{lemma}[theorem]{Lemma}

\newtheorem{proposition}[theorem]{Proposition}

\theoremstyle{definition}
\newtheorem{remark}[theorem]{Remark}
\newtheorem{definition}[theorem]{Definition}

\journal{}

\usepackage{geometry}

\begin{document}

\baselineskip15pt

\begin{frontmatter}



\title{\huge Interior $L_{p}$ regularity for Stokes systems  \tnoteref{t1}}
\tnotetext[t1]{Research supported by National Natural Science Foundation of China [grant 12001419 and grant 12071365], Fundamental Research Funds for the Central Universities [grant xzy012022004] and Shaanxi Fundamental Science Research Project for Mathematics and Physics
[grant 22JSZ003 and grant 23JSQ026].\\
\indent \mbox{}\hskip0.3cm All authors contributed to the study conception and design. \\
\indent \mbox{}\hskip0.3cm All authors declare that they have no potential conflict of interest.
}

\author[rvt1]{Rong Dong}
\ead{dongrong1203@mail.xjtu.edu.cn }
\author[rvt1]{Dongsheng Li\corref{cor1}}
\ead{lidsh@mail.xjtu.edu.cn} \cortext[cor1]{Corresponding author.}
\author[rvt2]{Lihe Wang}
\ead{lihe-wang@uiowa.edu}

\address[rvt1]{School of Mathematics and Statistics, Xi'an Jiaotong University, Xi'an 710049, China}
\address[rvt2]{Department of Mathematics, University of Iowa, Iowa City, IA52242, USA}

\begin{abstract}

A new iteration method is represented to study the interior $L_{p}$ regularity for Stokes systems both in divergence form and in non-divergence form. By the iteration, we improve the integrability of derivatives of solutions for Stokes systems step by step; after infinitely many steps, $L_{p}$ regularity is achieved; and in each step, the maximal function method is used where solutions and their derivatives are involved simultaneously in each scale. The H\"{o}lder continuity of the coefficients in spatial variables is assumed to compensate the different scalings between the solutions and their derivatives.

\end{abstract}

\bigskip


\end{frontmatter}


\section{Introduction}

We consider Stokes systems both in divergence form and in non-divergence form that are the following
\begin{equation}\label{s0.0}
\left\{
\begin{array}{ll}
\vspace{2mm} u_t+\mathcal{L}_1 u+\nabla \pi=\mbox{div}f~~\mbox{in}~~\Omega_{T},\\
\vspace{2mm}
\mbox{div}~u=g~~\mbox{in}~~\Omega_{T}
\end{array}
\right.
\end{equation}
and
\begin{equation}\label{s0.1}
\left\{
\begin{array}{ll}
\vspace{2mm} u_t+\mathcal{L}_2u+\nabla \pi=f~~\mbox{in}~~\Omega_{T},\\
\vspace{2mm}
\mbox{div}~u=g~~\mbox{in}~~\Omega_{T}
\end{array}
\right.
\end{equation}
respectively, with
\begin{equation*}
\mathcal{L}_1:=-\mbox{div}\Big(A(x,t)\nabla \Big)=:-D_i\Big(a^{ij}(x,t)D_j\Big)
\end{equation*}
and
\begin{equation*}
\mathcal{L}_2:=-A(x,t)D^2=:-a^{ij}(x,t)D_{ij},
\end{equation*}
where $1\leq i, j\leq n$ with $n\geq2$, $\Omega$ is a bounded domain in $\mathbb{R}^{n}$ and $\Omega_T=\Omega\times(-T,0]$ is a cylinder with $T>0$.
The summation convention is assumed throughout the paper.

Both Stokes systems in the form (\ref{s0.0}) and (\ref{s0.1}) arise from non-Newton fluid (c.f.  e.g. \cite{BP}) and also from Stokes systems on manifolds (c.f.  e.g. \cite{DM} and \cite{MT}).

Serrin's counterexample (cf. \cite{s14} and Remark 1.3 (iii) in the following) indicates that the interior regularity of Stokes systems is hard and will be essentially different from parabolic systems. The purpose of this paper is to investigate the interior $L_p$ regularity of them. For (\ref{s0.0}), we shall improve the integrability of $Du$, that is, we will start from $Du\in L_{p_0}$ and arrive at $Du\in L_{p}$ under the assumption that $f$, $g$ and $u\in L_{p}$ for some $1<p_0<p<\infty$. Here $u\in L_{p}$ is needed to provide the higher integrability in time direction, which can not be removed (cf. Remark 1.3 (iii) below) and this reveals that Stokes systems have very weak interior regularity in time direction. A new iteration method is designed so that the integrability of $Du$ is improved step by step and after infinitely many steps, the $L_p$ regularity (integrability) is achieved. In each step, the maximal function method is used and we gain the regularity from both parabolic equations and elliptic equations. Indeed, $\nabla\times u$ (\emph{see the notation list at the end of this section}) satisfies a parabolic equation and $Du$ can be estimated by $\nabla\times u$ and $u$ via an elliptic equation. Thus, $u$ and $\nabla\times u$ should be involved simultaneously, but unfortunately they have different scalings. To overcome this difficulty, we introduce a new and important quantity $\Theta$ which connects $u$ and $Du$ and whose scaling is between $u$ and $Du$. The H\"{o}lder continuity of the coefficients $a^{ij}$ in the spatial variables is assumed to compensate the different scalings of these quantities. For (1.2), we will improve the integrability of $D^2u$ from $L_{p_0}$ to $L_{p}$ for some $1<p_0<p<\infty$ by the same idea.

The interior $L_p$ regularity of Stokes systems has not been well developed yet.
As for constant coefficients case, Hu, Li and Wang \cite{s11} establish interior $L_p$ regularity for $\Delta u$, that is, $\Delta u\in L_p$ is attained under the assumption $u\in W^{2,1}_2$ and righthand side $f\in L_p$ with $p>2$; and if $u\in L_p$ is further assumed, $D^2 u\in L_p$ can be obtained.
Observe that in \cite{s11} since the coefficients of the Stokes system are constant already, we do not need to transform it to a new system with coefficients depending only on $t$ and then the righthand side having the derivatives of $u$ as what is done in this paper. Consequently in \cite{s11}, while using maximal function method, in each scale, only $\Delta u$ and $f$ are involved which have the same scaling. Therefore, the iteration argument designed in this paper is not needed there. In \cite{s24}, Meyers improved integrability of the first derivatives of solutions of elliptic equations in divergence form from $L_2$ to $L_p$ for some $p>2$ which may be very closed to 2, where the righthand sides of the equations are in the space $W^{-1}_p$ and no continuity of the coefficients is needed.
By this technique,
Kaplick$\acute{y}$ and Wolf \cite{s12} generalize Meyers' result to Stokes systems in divergence form and they actually show that under the assumption $u\in L_{\infty}(0,T;L_{2}(\Omega))\cap L_{2}(0,T;W^{1}_{2}(\Omega))$, there exists some
$p\in[2,q]$ with $2<q<\frac{2(n+2)}{n}$ such that $Du\in L_p$. We should mention here that under their assumption,
one can deduce $u\in L_{2(n+2)/{n}}(\Omega_T)$ by Sobolev embedding theorem and this implies that the upper bound of $q$ is $\frac{2(n+2)}{n}$.

For a priori $L_p$ estimates, Solonnikov \cite{s15} is the first important contribution which demonstrates a global a priori $L_p$ estimate for Stokes systems with constant coefficients by using potential theory. Later, Solonnikov \cite{s21} further develops the result to the case for variable coefficients that are continuous and belong to $W_q^1$ in the spatial variables for some large $q$. Recently, Dong and Phan \cite{s9} obtain an a priori interior mixed-norm estimate via sharp function method, where the coefficients are assumed to have small mean oscillations in the spatial variables. In these a priori estimates, $Du$ and $D^2u\in L_p$ are assumed for the systems in divergence form and in non-divergence form respectively, and the assumptions are essentially used in the proofs.

\bigskip

In this paper, we make the following assumptions for the coefficient matrix $A(y,s)=\big(a^{ij}(y,s)\big)$ :

\bigskip

($\bf H_1$) Ellipticity:

There exists $\lambda\in(0,1]$ such that for any $(y,s)\in\Omega_T$ and $\xi=(\xi_i)\in\mathbb{R}^{n}$,
\begin{equation*}
a^{ij}(y,s)\xi_i\xi_j\geq\lambda|\xi|^2, ~~|a^{ij}(y,s)|\leq\lambda^{-1}.
\end{equation*}

($\bf H_2$) Regularity:

There exist $B>0$ and $0<\gamma<1$ such that
\begin{equation*}
[A]_{C_{x}^{\gamma}(\Omega_T)}\leq B.
\end{equation*}

Here in the ellipticity condition ($H_1$), we do not need $A$ is symmetric. In the regularity condition ($H_2$), we only need the H\"{o}lder continuity of $A$ in space direction (\emph{please refer to the notation list at the end of this section}). One should notice that the regularity of $u_t$ can not be obtained even if $A$ is smooth in all variables.

\vspace{2mm}
Now we state our interior $L_{p}$ regularity corresponding to (\ref{s0.0}) and (\ref{s0.1}) respectively.

\begin{theorem}
Let $1<p_0< p<\infty$. Suppose $u\in W_{p_0}^{1,\frac12}(\Omega_{T};\mathbb{R}^n)\cap L_p(\Omega_{T};\mathbb{R}^n)$, $\pi\in L_{1}(\Omega_{T})$, $f\in L_{p}(\Omega_{T};\mathbb{R}^{n^2})$ and $g\in L_{p}(\Omega_{T})$ satisfying (\ref{s0.0}) in weak sense. Assume that $(H_1)$ and $(H_2)$ hold.
Then we have $Du\in L_{p,loc}(\Omega_{T};\mathbb{R}^{n^2})$ and for any $\tilde\Omega_{T}=\tilde\Omega\times(a,b]$ with $\tilde\Omega\subset\subset\Omega$ and $-T<a<b\leq0$,
\begin{eqnarray}\label{s.01}
\|Du\|_{L_p(\tilde\Omega_{T})}
\leq
C\Big(\|u\|_{L_p(\Omega_{T})}+\|f\|_{L_p(\Omega_{T})}+\|g\|_{L_p(\Omega_{T})}\Big),
\end{eqnarray}
where $C$ depends only on $n$, $\lambda$, $p_0$, $p$, $B$, $\gamma$, $\tilde\Omega_{T}$ and $\Omega_{T}$.
\end{theorem}

\begin{theorem}
Let $1<p_0<p<\infty$. Suppose $u\in W_{p_0}^{2,1}(\Omega_{T};\mathbb{R}^n)\cap L_p(\Omega_{T};\mathbb{R}^n)$, $\nabla \pi\in L_{1}(\Omega_{T};\mathbb{R}^n)$, $f\in L_{p}(\Omega_{T};\mathbb{R}^{n})$ and $Dg\in L_{p}(\Omega_{T};\mathbb{R}^{n})$ satisfying (\ref{s0.1}). Assume that $(H_1)$ and $(H_2)$ hold.
Then we have $D^2u\in L_{p,loc}(\Omega_{T};\mathbb{R}^{n^3})$ and for any $\tilde\Omega_{T}=\tilde\Omega\times(a,b]$ with $\tilde\Omega\subset\subset\Omega$ and $-T<a<b\leq0$,
\begin{eqnarray}\label{s.02}
\|D^2u\|_{L_p(\tilde\Omega_{T})}
\leq
C\Big(\|u\|_{L_p(\Omega_{T})}
+\|f\|_{L^p(\Omega_{T})}+\|Dg\|_{L_p(\Omega_{T})}\Big),
\end{eqnarray}
where $C$ depends only on $n$, $\lambda$, $p_0$, $p$, $B$, $\gamma$, $\tilde\Omega_{T}$ and $\Omega_{T}$.
\end{theorem}

\begin{remark}
(i) In Theorem 1.1, since $u\in W_{p_0}^{1,\frac12}(\Omega_{T})$ implies $u\in L_{\frac{(n+2)p_0}{n+2-p_0}}(\Omega_{T})$ by Sobolev embedding, we have (\ref{s.01}) holds with $p$ replaced by $\frac{(n+2)p_0}{n+2-p_0}$ without the assumption $u\in L_p(\Omega_{T})$. Similarly, in Theorem 1.2, we have (\ref{s.02}) holds with $p$ replaced by $\frac{(n+2)p_0}{n+2-2p_0}$ without the assumption $u\in L_p(\Omega_{T})$. Even for these cases, they are still new regularity results.

(ii) If we strengthen the assumption to $u\in W_p^{1,\frac12}(\Omega_T)$ in Theorem 1.1 or to $u\in W_p^{2,1}(\Omega_T)$ in Theorem 1.2, the Condition ($H_2$) can be weakened to
small mean oscillations with respect to the spatial variable
to obtain (\ref{s.01}) or (\ref{s.02}) (cf. \cite{s9}). Nevertheless, those are actually a priori estimates.

(iii) Let $h$ be a harmonic function in $\Omega$ and $\sigma=(\varepsilon_0-t)^{-\frac{\varepsilon_0}{p_0}}$ with $-1<t\leq 0$, $0<\varepsilon_0<1$ and $1<p_0<p$. Then $u(x,t)=\nabla h(x)\sigma(t)$ and $\pi(x,t)=-h(x)\sigma'(t)$ is a pair of solution of (1.1) or (1.2) with $a^{ij}=\delta^{ij}$ and vanishing $f$ and $g$ (cf. \cite{s14} and \cite{s11}). It is clear that $\int_{-1}^0\sigma^{p_0}(t)dt<2$ and $\int_{-1}^0\sigma'(t)dt<2$, that is, the pressure has a uniform bound in $L^1$ norm but $\int_{-1}^0\sigma^{p}(t)dt$ has no uniform bound as $p=\frac{2p_0}{\varepsilon_0}$ and $\varepsilon_0$ small. It follows that $\|Du\|_{L_p(\tilde{\Omega}\times(-\frac{1}{2},0])}$ and $\|D^2u\|_{L_p(\tilde{\Omega}\times(-\frac{1}{2},0])}$ have no uniform bound. Then the term $\|u\|_{L_p(\Omega_T)}$ can not be removed or weakened in both (\ref{s.01}) and (\ref{s.02}). This illustrates that Stokes systems inherently have weak interior regularity in time direction. The $L_p$ integrability of $u$ in both (\ref{s.01}) and (\ref{s.02}) is for obtaining the integrability of $Du$ and $D^2u$ in time direction. Also $u_t$ has no more integrability even if $A$ is smooth.

\end{remark}

This paper is organized as follows. In Section 2, we  introduce some important tools. Section 3 and 4 are devoted to the proof of Theorem 1.1 and 1.2 respectively, where our new iteration is demonstrated.

\vspace{2mm}

Before ending this section, we introduce some notations.
\vskip 2mm\noindent{\bf Notations:}
{\baselineskip=20pt

$B_r(x):=\{y\in\mathbb{R}^n:|x-y|<r\}$;
~$B_r:=B_r(0)$.

$Q_r(x,t):=B_r(x)\times(t-r^2,t]$;~
$Q_r:=B_r\times(-r^2,0]$.

$K_r(x,t):=B_r(x)\times(t-\frac{r^2}{2},t+\frac{r^2}{2}]$;~
$K_r:=B_r\times(-\frac{r^2}{2},\frac{r^2}{2}]$.

$\overline{u}_{Q_r(x,t)}:=\fint_{Q_r(x,t)}u(y,s)
dyds=\frac{1}{|Q_r(x,t)|}\int_{Q_r(x,t)}u(y,s)
dyds$.

$\overline{u}_{B_r(x)}(t):=\fint_{B_r(x)}u(y,t)
dy=\frac{1}{|B_r(x)|}\int_{B_r(x)}u(y,t)
dy$.

$D_iu:=\frac{\partial u}{\partial x_i}$; $Du:=\big(D_{1}u,D_{2}u,...,D_{n}u\big)$; $D^2u:=\big(D_iD_ju\big)_{n\times n}$.

$\nabla\times u:=\big(D_ju_i-D_iu_j\big)_{n\times n}$.

$[u]_{C_{x}^{\gamma}(\Omega_{T})}:=\sup\limits_{(x,t),(y,t)\in\Omega_{T}, x\neq y}\frac{|u(x,t)-u(y,t)|}{|x-y|^{\gamma}}$.

$\partial_p \Omega_T: \mbox{parabolic~boundary~of}~\Omega_T$.
}

\section{Preliminary}
In this section, we first introduce some analytic tools in parabolic spaces, which can be found in \cite{s2}.

\begin{lemma}
Suppose $f$ is a measurable function in a bounded domain $\Omega_T\subset\mathbb{R}^{n+1}$. Let $\theta>0$ and $m>1$ be constants. Then for $0<p<\infty$,
\begin{eqnarray*}
f\in L_p(\Omega_T) ~\Leftrightarrow~ S:=\sum_{k\geq 1}m^{kp}|\{(x,t)\in\Omega_T: |f(x,t)|>\theta m^k\}|<\infty
\end{eqnarray*}
and
\begin{eqnarray*}
\frac{1}{C}S\leq\|f\|^p_{L_p(\Omega_T)}\leq C(|\Omega_T|+S),
\end{eqnarray*}
where $C$ depends only on $\theta$, $m$ and $p$.
\end{lemma}

\begin{definition}
Let $f$ be an integrable function on $\Omega_T$. The Hardy-Littlewood maximal function is defined by
\begin{eqnarray*}
\mathcal{M}_{\Omega_T}(f)(x,t)=\sup_{r>0}\fint_{K_r(x,t)\cap\Omega_T}|f(y,s)|dyds.
\end{eqnarray*}
We will drop the index $\Omega_T$ if $\Omega_T$ is understood clearly in the context.
\end{definition}

\begin{lemma}
(i) If $f\in L_p(\Omega_T)$ with $p>1$, then $\mathcal{M}(f)\in L_p(\Omega_T)$. Moreover,
\begin{eqnarray*}
\|\mathcal{M}(f)\|_{L_p(\Omega_T)}\leq C\|f\|_{L_p(\Omega_T)},
\end{eqnarray*}
where $C$ depends only on $n$ and $p$;

(ii) If $f\in L_1(\Omega_T)$, then for any $\xi>0$,
\begin{eqnarray*}
|\{(x,t)\in\Omega_T:\mathcal{M}(f)(x,t)>\xi\}|\leq\frac{C}{\xi}\|f\|_{L_1(\Omega_T)},
\end{eqnarray*}
where $C$ depends only on $n$.
\end{lemma}

We also need the following modified Vitali covering lemma (cf. \cite{s2}), which is crawling of ink spots lemma in \cite{s26}.

\begin{lemma}
Let $0<\varepsilon<1$ and $E\subset F\subset Q_1$ be two measurable sets such that
\begin{eqnarray*}
|E|<\varepsilon|Q_1|,
\end{eqnarray*}
and for any $(x,t)\in Q_1$ and $r\in (0,1]$, $|E\cap K_r(x,t)|\geq\varepsilon|K_r(x,t)|$ implies $$K_r(x,t)\cap Q_1\subset F.$$
Then
\begin{eqnarray*}
|E|\leq10^{n+2}\varepsilon|F|.
\end{eqnarray*}
\end{lemma}

\bigskip

Next, we give the following $L_p$ regularity for parabolic equations with bounded coefficients depending only on $t$. Some basic results for discontinuous coefficients can be found in \cite{s23}.

\begin{lemma}
Let $1<p<\infty$ and $\tilde\Omega_{T}\subset\subset\Omega_T$. Suppose $(a^{ij}(t))$ satisfies $(H_1)$, $u\in L_p(-T,0;W^{-1}_p(\Omega))$, $u_t\in L_p(-T,0;W^{-3}_p(\Omega))$ and $g\in L_p(-T,0;W^{-3}_p(\Omega))$ such that
\begin{equation*}
u_t-a^{ij}(t)D_{ij}u=g~~\mbox{in}~~L_p(-T,0;W^{-3}_p(\Omega)).
\footnote{Let $<,>$ be the pairing between $W^{-3}_p(\Omega)$ and its dual $W^{3}_{p',0}(\Omega)$
with $\frac1p+\frac{1}{p'}=1$. The equation means $\int^0_{-T} <u_t-a^{ij}(t)D_{ij}u,\varphi>dt=\int^0_{-T} <g,\varphi> dt$ for any $\varphi\in L_{p'}(-T,0;W^{3}_{p',0}(\Omega))$.}
\end{equation*}

(i) If $u\in L_p(\Omega_{T})$ and $g=0$, then $u,Du\in L_{\infty,loc}(\Omega_T)$ and
\begin{eqnarray*}
\|u\|_{L_\infty(\tilde{\Omega}_T)}+\|Du\|_{L_\infty(\tilde{\Omega}_T)}\leq C\|u\|_{L_p(\Omega_{T})},
\end{eqnarray*}
where $C$ depends only on $n$, $\lambda$, $p$, $\tilde\Omega_{T}$ and $\Omega_{T}$.

(ii) If $g=div f$ for some $f\in L_{p}(\Omega_{T};\mathbb{R}^{n})$, then $Du\in L_{p,loc}(\Omega_T)$ and
\begin{eqnarray}\label{r.2}
\|Du\|_{L_p(\tilde{\Omega}_T)}\leq C\Bigg\{\|f\|_{L_p(\Omega_{T})}+\Big(\int_{-T}^0\| u(\cdot,t)\|^{p}_{W^{-1}_p(\Omega)}dt\Big)^{\frac{1}{p}}\Bigg\};
\end{eqnarray}
and moreover, if $u=0~\mbox{on}~\partial_p\Omega_{T}$, then
\begin{eqnarray}\label{r.4}
\|Du\|_{L_p(\tilde{\Omega}_T)}\leq C\|f\|_{L_p(\Omega_{T})},
\end{eqnarray}
where $C$ depends only on $n$, $\lambda$, $p$, $\tilde\Omega_{T}$ and $\Omega_{T}$.

(iii) If $g=D_{ij}f_{ij}$ for some $f=(f_{ij})\in L_{p}(\Omega_{T};\mathbb{R}^{n^2})$, then $u\in L_{p,loc}(\Omega_T)$ and
\begin{eqnarray}\label{r.5}
\|u\|_{L_p(\tilde{\Omega}_T)}\leq C\Bigg\{\|f\|_{L_p(\Omega_{T})}+\Big(\int_{-T}^0\| u(\cdot,t)\|^{p}_{W^{-1}_p(\Omega)}dt\Big)^{\frac{1}{p}}\Bigg\};
\end{eqnarray}
and moreover, if $u=0~\mbox{on}~\partial_p\Omega_{T}$\footnote{Here $u$ satisfies the equation with vanishing boundary condition means
$$\int_{\Omega_{T}} u \big(-\varphi_t-a^{ij}(t)D_{ij}\varphi\big)dxdt=\int_{\Omega_{T}} f_{ij}D_{ij}\varphi dxdt$$ for any $\varphi\in{\cal{T}}=\{\varphi\in C^2(\Omega_{T})\cap C(\overline{\Omega_{T}}):\varphi(x,0)=0~\mbox{for}~x\in\Omega~\mbox{and}~ \varphi(x,t)=0~ \mbox{for}~ (x,t)\in\partial\Omega\times[-T,0]\}$. },
then
\begin{eqnarray}\label{r.3}
\|u\|_{L_p(\tilde{\Omega}_T)}\leq C\|f\|_{L_p(\Omega_{T})},
\end{eqnarray}
where $C$ depends only on $n$, $\lambda$, $p$, $\tilde\Omega_{T}$ and $\Omega_{T}$.
\end{lemma}

\section{Regularity for Stokes systems in divergence form}\setcounter{equation}{0}
In this section, we will show Theorem 1.1, the interior $L_{p}$ regularity for Stokes systems in divergence form. By suitable dilation, instead of $(H_2)$, we assume $(H_2')$:
\begin{equation*}
[A]_{C_{x}^{\gamma}(\Omega_T)}\leq \delta.
\end{equation*}
where $0<\delta<1$ and will be determined in the following. By standard scaling and covering arguments, it suffices to prove:

\begin{theorem}
Let $1<p_0< p<\infty$. Suppose $u\in W_{p_0}^{1,\frac12}(\Omega_{T};\mathbb{R}^n)\cap L_p(\Omega_{T};\mathbb{R}^n)$, $\pi\in L_{1}(\Omega_{T})$, $f\in L_{p}(\Omega_{T};\mathbb{R}^{n^2})$ and $g\in L_{p}(\Omega_{T})$ satisfying (\ref{s0.0}) in weak sense. Assume that $Q_{7}\subset\Omega_T$, $(H_1)$ and $(H_2')$ hold with some $0<\delta<1$ depending only on $n$, $\lambda$, $\gamma$, $p_0$ and $p$.
Then we have $Du\in L_{p}(Q_{\frac{1}{2}};\mathbb{R}^{n^2})$ and
\begin{eqnarray}\label{n1}
\|Du\|_{L_p(Q_{\frac12})}
\leq
C\Big(\|u\|_{L_p(Q_7)}+\|f\|_{L_p(Q_7)}+\|g\|_{L_p(Q_7)}\Big),
\end{eqnarray}
where $C$ depends only on $n$, $\lambda$, $\gamma$, $p_0$ and $p$.
\end{theorem}

\bigskip

As mentioned before, we will use an iteration to prove Theorem 3.1. To establish the iteration inequality (see (\ref{d10.49})), the following proposition is the key, in whose proof regularity of both parabolic equations and elliptic equations will be applied. From now on, we shall abbreviate $\mathcal{M}_{Q_7}(f)=\mathcal{M}(f)$ for any function $f$ in this and the next section.

\bigskip

\begin{proposition}
Set $0<\varepsilon<1<p_0<\infty$ and $q_0=\frac{1+p_0}{2}$. Suppose
$u\in W_{p_0}^{1,\frac{1}{2}}(\Omega_{T};\mathbb{R}^{n})$, $\pi\in L_{1}(\Omega_{T})$, $f\in L_{p_0}(\Omega_{T};\mathbb{R}^{n^2})$, $g\in L_{p_0}(\Omega_{T})$ satisfying (\ref{s0.0}) in weak sense, $Q_7\subset\Omega_{T}$ and $(H_1)$ holds. There exists a constant $0<\tilde{\delta}<1$ depending only on $\varepsilon$, $n$, $\lambda$ and $p_0$. Assume that $(H_2')$ holds with $0<\delta\leq\tilde{\delta}$. For any $(x,t)\in Q_1$, let
\begin{eqnarray}\label{d4.37}
\Theta(x,t):=\sup_{r>0}r^{q_0(\gamma-1)}\fint_{K_r(x,t)\cap Q_7}|u(y,s)-\overline{u}_{B_r(x)}(s)|^{q_0}dyds,
\end{eqnarray}
where $\gamma$ is given in ($H_2'$).
There exists a constant $N_1>1$ depending only on $n$, $\lambda$ and $p_0$ such that if
\begin{eqnarray*}
|\{(x,t)\in Q_1:\mathcal{M}(|\nabla\times u|^{q_0})>N_1^{q_0}\}|\leq \varepsilon|Q_1|,
\end{eqnarray*}
then
\begin{equation*}
\begin{aligned}
\lefteqn{|\{(x,t)\in Q_1:\mathcal{M}(|\nabla\times u|^{q_0})>N_1^{q_0}\}|}\hspace*{16mm}\\
\leq&\ 10^{n+2}\varepsilon\Big(|\{(x,t)\in Q_1:\mathcal{M}(|\nabla\times u|^{q_0})>1\}|+|\{(x,t)\in Q_1:\Theta>1\}|\\
& +|\{(x,t)\in Q_1:\mathcal{M}(|f|^{q_0})>\tilde{\delta}^{q_0}\}|+|\{(x,t)\in Q_1:\mathcal{M}(|g|^{q_0})>1\}|\Big).
\end{aligned}
\end{equation*}
\end{proposition}

\begin{remark}
The scaling of $\Theta$ is between $u$ and $Du$. In the proof of Proposition 3.2, $\Theta$ will be used to control $u$. In Lemma 3.6, $\Theta$ will be estimated by $u$ and $Du$.
\end{remark}

\bigskip

Lemma 2.4 will be utilized to prove Proposition 3.2. The following lemma is the key
to establish the second condition in Lemma 2.4.

\bigskip

\begin{lemma}
Under the assumptions of Proposition 3.2, there exists a constant $N_1>1$ depending only on $n$, $\lambda$ and $p_0$ such that if
\begin{eqnarray}\label{dl.20}
\begin{aligned}
\lefteqn{\{(x,t)\in Q_1:\mathcal{M}(|\nabla\times u|^{q_0})\leq1\}\cap\{(x,t)\in Q_1:\Theta\leq1\}}\hspace*{3mm}\\
& \cap\{(x,t)\in Q_1:\mathcal{M}(|f|^{q_0})\leq\tilde{\delta}^{q_0}\}\cap\{(x,t)\in Q_1:\mathcal{M}(|g|^{q_0})\leq1\}\cap K_r(\hat{x},\hat{t})
\neq\emptyset
\end{aligned}
\end{eqnarray}
for some $(\hat{x},\hat{t})\in Q_1$ and $0<r\leq1$, then
\begin{eqnarray}\label{d8.20}
|\{(x,t)\in Q_1:\mathcal{M}(|\nabla\times u|^{q_0})>N^{q_0}_1\}\cap K_r(\hat{x},\hat{t})|\leq\varepsilon|K_r(\hat{x},\hat{t})|.
\end{eqnarray}
\end{lemma}

\vskip 2mm\noindent{\bf Proof.}\quad
Let $K_r(\hat{x},\hat{t})$ be given such that (\ref{dl.20}) holds. Since $\Theta$ has different scaling from $\nabla\times u$, $f$ and $g$, we can not scale $K_r(\hat{x},\hat{t})$ into $K_1$ as usual.

In view of (\ref{dl.20}), there exists a point $(\bar{x},\bar{t})\in Q_{1}\cap K_r(\hat{x},\hat{t})$ so that for any $\rho>0$,
\begin{eqnarray}\label{dl.21}
\fint_{K_{\rho}(\bar{x},\bar{t})\cap Q_7}|\nabla\times u|^{q_0}dxdt\leq1,
\end{eqnarray}
\begin{eqnarray}\label{dl.22}
\fint_{K_{\rho}(\bar{x},\bar{t})\cap Q_7}|f|^{q_0}dxdt\leq\tilde{\delta}^{q_0},
~~\fint_{K_{\rho}(\bar{x},\bar{t})\cap Q_7}|g|^{q_0}dxdt\leq1
\end{eqnarray}
and (by $\Theta\leq1$)
\begin{eqnarray}\label{dl.27}
\fint_{K_{\rho}(\bar{x},\bar{t})\cap Q_7}|u-\overline{u}_{B_{\rho}(\bar{x})}|^{q_0}dxdt\leq\rho^{q_0(1-\gamma)}.
\end{eqnarray}

In the following, we abbreviate $Q_\rho(\hat{x},\tilde{t})=\widetilde{Q}_{\rho}$ with $\tilde{t}=\min\{\hat{t}+3r^2,0\}$ for any $\rho>0$. By $\widetilde{Q}_{5r}\subset K_{8r}(\bar{x},\bar{t})\cap Q_7$ and (\ref{dl.21})-(\ref{dl.27}),
\begin{eqnarray}\label{dl.26}
\fint_{\widetilde{Q}_{5r}}|\nabla\times u|^{q_0}dxdt\leq C,
\end{eqnarray}
\begin{eqnarray}\label{dl.28}
\fint_{\widetilde{Q}_{5r}}|f|^{q_0}dxdt\leq C\tilde{\delta}^{q_0},
~~\fint_{\widetilde{Q}_{5r}}|g|^{q_0}dxdt\leq C
\end{eqnarray}
and
\begin{equation}\label{dl.66}
\begin{aligned}
\fint_{\widetilde{Q}_{5r}}|u-\overline{u}_{B_{5r}(\hat{x})}|^{q_0}dxdt
\leq Cr^{q_0(1-\gamma)},
\end{aligned}
\end{equation}
where $C$ depends only on $n$.

Recall $\overline{A}_{B_r(x)}(t)=\fint_{B_r(x)}A(y,t)dy$. Take $\nabla\times$ on both sides of the first equation in (\ref{s0.0}) and then, by $\nabla\times\nabla \pi=0$,
\begin{equation}\label{n.78}
(\nabla\times u)_t-\mbox{div}\big(\overline{A}_{B_{4r}(\hat{x})}(t)\nabla(\nabla\times u)\big)
=\nabla\times\mbox{div}~\tilde{f}~~\mbox{in}~~\widetilde{Q}_{4r},\\
\end{equation}
where $\tilde{f}=\big(A-\overline{A}_{B_{4r}(\hat{x})}(t)\big)Du+f$.
Let $u_\varepsilon$ and $\tilde{f}_\varepsilon$ be the mollification of $u$ and $\tilde{f}$ in $x$. Then
\begin{equation*}
(\nabla\times u_\varepsilon)_t-\mbox{div}\big(\overline{A}_{B_{4r}(\hat{x})}(t)\nabla(\nabla\times u_\varepsilon)\big)
=\nabla\times\mbox{div}~\tilde{f}_\varepsilon~~\mbox{in}~~\widetilde{Q}_{4r}.
\end{equation*}
We separate $\nabla\times u_\varepsilon$ into $w_\varepsilon+v_\varepsilon$ with $w_\varepsilon$ solving
\begin{equation*}
\left\{
\begin{array}{ll}
(w_\varepsilon)_t-\mbox{div}\big(\overline{A}_{B_{4r}(\hat{x})}(t)\nabla w_\varepsilon\big)
=\nabla\times \mbox{div}~\tilde{f}_\varepsilon~~\mbox{in}~~\widetilde{Q}_{4r},\\
\vspace{2mm}
w_\varepsilon=0~~\mbox{on}~~\partial_p\widetilde{Q}_{4r}.
\end{array}
\right.
\end{equation*}
It follows that
\begin{equation*}
(v_\varepsilon)_t-\mbox{div}\big(\overline{A}_{B_{4r}(\hat{x})}(t)\nabla v_\varepsilon\big)
=0~~\mbox{in}~~\widetilde{Q}_{4r}.
\end{equation*}

Let $\lim\limits_{\varepsilon\rightarrow0^+} w_\varepsilon=w$ and $\lim\limits_{\varepsilon\rightarrow0^+} v_\varepsilon=v$. Then $\nabla\times u=w+v$. Using (\ref{r.3}) to the equation that $w_\varepsilon$ satisfies, we obtain
\begin{equation}\label{dl.37}
\begin{aligned}
\fint_{\widetilde{Q}_{\frac{7r}{2}}}|w|^{q_0}dxdt
=&\ \lim_{\varepsilon\rightarrow0^+}\fint_{\widetilde{Q}_{\frac{7r}{2}}}|w_\varepsilon|^{q_0}dxdt
\leq C\lim_{\varepsilon\rightarrow0^+}\fint_{\widetilde{Q}_{4r}}|\tilde{f}_\varepsilon|^{q_0}dxdt\\
\leq&\ C\fint_{\widetilde{Q}_{4r}}|\big(A-\overline{A}_{B_{4r}(\hat{x})}(t)\big)D u|^{q_0}dxdt+C\fint_{\widetilde{Q}_{4r}}|f|^{q_0}dxdt,
\end{aligned}
\end{equation}
where $C$ depends only on $n$, $\lambda$ and $p_0$.
Letting $\varepsilon\rightarrow0^+$, we also obtain
\begin{equation*}
v_t-\mbox{div}\big(\overline{A}_{B_{4r}(\hat{x})}(t)\nabla v\big)
=0~~\mbox{in}~~\widetilde{Q}_{4r}.
\end{equation*}
Since,
\begin{eqnarray}\label{n.77}
\begin{aligned}
\Delta u
=\mbox{div}~\nabla\times u+\nabla~\mbox{div}~u,
\end{aligned}
\end{eqnarray}
we have, by $\mbox{div}~u=g$,
\begin{eqnarray*}
\begin{aligned}
\Delta \big(u-\overline{u}_{B_{5r}(\hat{x})}\big)
=\mbox{div}~\nabla\times u+\nabla g.
\end{aligned}
\end{eqnarray*}
Using $W^{1,p}$ estimates to $u-\overline{u}_{B_{5r}(\hat{x})}$
and by (\ref{dl.26})-(\ref{dl.66}),
\begin{eqnarray*}
\begin{aligned}
\fint_{\widetilde{Q}_{4r}}|Du|^{q_0}dxdt
=&\ \fint_{\widetilde{Q}_{4r}}|D(u-\overline{u}_{B_{5r}(\hat{x})})|^{q_0}dxdt\\
\leq&\ C\Big(\fint_{\widetilde{Q}_{5r}}|\nabla\times u|^{q_0}dxdt+\fint_{\widetilde{Q}_{5r}}| g|^{q_0}dxdt\\
&\mbox{}\hskip1.5cm +r^{-q_0}\fint_{\widetilde{Q}_{5r}}|u-\overline{u}_{B_{5r}(\hat{x})}|^{q_0}dxdt\Big)
\leq C\Big(1+r^{{-q_0}\gamma}\Big),
\end{aligned}
\end{eqnarray*}
where $C$ depends only on $n$ and $p_0$.
Observe the righthand side is unbounded as $r\rightarrow0$. However, we can use the H\"{o}lder continuity of the coefficient $A$ to offset the unboundedness. Actually, by $(H_2')$, it follows that
\begin{eqnarray*}
\begin{aligned}
\fint_{\widetilde{Q}_{4r}}|\big(A-\overline{A}_{B_{4r}(\hat{x})}(t)\big)Du|^{q_0}dxdt
\leq&\ \|A-\overline{A}_{B_{4r}(\hat{x})}(t)\|^{q_0}_{L_\infty(\widetilde{Q}_{4r})}\fint_{\widetilde{Q}_{4r}}|Du|^{q_0}dxdt\\
\leq&\ Cr^{{q_0}\gamma}\delta^{q_0}\Big(1+r^{{-q_0}\gamma}\Big)
\leq C\tilde{\delta}^{q_0}.
\end{aligned}
\end{eqnarray*}
Combining the above inequality with (\ref{dl.28}) and (\ref{dl.37}), we obtain
\begin{eqnarray}\label{dl.29}
\fint_{\widetilde{Q}_{\frac{7r}{2}}}|w|^{q_0}dxdt\leq C\tilde{\delta}^{q_0},
\end{eqnarray}
where $C$ depends only on $n$, $\lambda$ and $p_0$.

Since, by (\ref{dl.26}) and (\ref{dl.29}),
\begin{eqnarray*}
\begin{aligned}
\fint_{\widetilde{Q}_{\frac{7r}{2}}}|v|^{q_0}dxdt
= \fint_{\widetilde{Q}_{\frac{7r}{2}}}|\nabla\times u-w|^{q_0}dxdt
\leq C\fint_{\widetilde{Q}_{\frac{7r}{2}}}\Big(|\nabla\times u|^{q_0}+|w|^{q_0}\Big)dxdt
\leq C,
\end{aligned}
\end{eqnarray*}
we have, by Lemma 2.5 (i),
\begin{eqnarray}\label{dl.31}
\|v\|^{q_0}_{L_\infty(\widetilde{Q}_{3r})}\leq C\fint_{\widetilde{Q}_{\frac{7r}{2}}}|v|^{q_0}dxdt\leq N^{q_0}_0,
\end{eqnarray}
where $N_0$ is a constant depending only on $n$, $\lambda$ and $p_0$.

Let $N^{q_0}_1=\max\{2^{q_0}N^{q_0}_0,4^{n+2}\}$.
We claim that
for any $(x,t)\in Q_1\cap K_r(\hat{x},\hat{t})$,
\begin{eqnarray}\label{dl.33}
\mathcal{M}_{\widetilde{Q}_{\frac{7r}{2}}}(|w|^{q_0})(x,t)\leq N_0^{q_0}\Longrightarrow
\mathcal{M}(|\nabla\times u|^{q_0})(x,t)\leq N^{q_0}_1.
\end{eqnarray}
Actually, suppose $\mathcal{M}_{\widetilde{Q}_{\frac{7r}{2}}}(|w|^{q_0})(x,t)\leq N_0^{q_0}$ and we show $\mathcal{M}(|\nabla\times u|^{q_0})(x,t)\leq N^{q_0}_1$ by two cases.
If $0<\rho\leq2r$, then $K_{\rho}(x,t)\cap Q_7\subset \widetilde{Q}_{3r}$. By $\nabla\times u=w+v$ and (\ref{dl.31}), we have
\begin{eqnarray*}
\fint_{K_{\rho}(x,t)\cap Q_7}|\nabla\times u|^{q_0}dxdt
\leq2^{q_0-1}\fint_{K_{\rho}(x,t)\cap Q_7}\Big(|w|^{q_0}+|v|^{q_0}\Big)dxdt\leq 2^{q_0}N^{q_0}_0\leq N^{q_0}_1.
\end{eqnarray*}
If $\rho>2r$, then $K_{\rho}(x,t)\cap Q_7\subset K_{2\rho}(\bar{x},\bar{t})\cap Q_7$.
Clearly, (\ref{dl.21}) implies
\begin{eqnarray*}
\begin{aligned}
\fint_{K_{\rho}(x,t)\cap Q_7}|\nabla\times u|^{q_0}dxdt
\leq 4^{n+2}\fint_{K_{2\rho}(\bar{x},\bar{t})\cap Q_7}|\nabla\times u|^{q_0}dxdt
\leq 4^{n+2}\leq N^{q_0}_1.
\end{aligned}
\end{eqnarray*}
Hence, in view of Definition 2.2, (\ref{dl.33}) is true and it implies
\begin{eqnarray}\label{n.63}
\begin{aligned}
\lefteqn{\{(x,t)\in Q_1:\mathcal{M}(|\nabla\times u|^{q_0})>N^{q_0}_1\}\cap K_r(\hat{x},\hat{t})}\hspace*{12mm}\\
&\subset\{(x,t)\in Q_1:\mathcal{M}_{\widetilde{Q}_{\frac{7r}{2}}}(|w|^{q_0})>N^{q_0}_0\}\cap K_r(\hat{x},\hat{t}).
\end{aligned}
\end{eqnarray}

From (\ref{n.63}), Lemma 2.3 (ii) and (\ref{dl.29}), we get
\begin{eqnarray*}
\begin{aligned}
\lefteqn{|\{(x,t)\in Q_1:\mathcal{M}(|\nabla\times u|^{q_0})>N^{q_0}_1\}\cap K_r(\hat{x},\hat{t})|}\hspace*{12mm}\\
&\leq\ |\{(x,t)\in Q_1:\mathcal{M}_{\widetilde{Q}_{\frac{7r}{2}}}(|w|^{q_0})>N_0^{q_0}\}\cap K_r(\hat{x},\hat{t})|\\
&\mbox{}\hskip1cm\leq\ \frac{C}{N^{q_0}_0}\int_{\widetilde{Q}_{\frac{7r}{2}}}|w|^{q_0}dxdt
\leq \frac{C}{N^{q_0}_0}\tilde{\delta}^{q_0}|\widetilde{Q}_{\frac{7r}{2}}|,
\end{aligned}
\end{eqnarray*}
where $C$ depends only on $n$, $\lambda$ and $p_0$.
Taking $\tilde{\delta}$ small enough, this leads to (\ref{d8.20}) clearly.
\qed

\bigskip

\noindent{\bf Proof of Proposition 3.2.}\quad
Since Lemma 3.4 is equivalent to that if
\begin{eqnarray*}
|\{(x,t)\in Q_1:\mathcal{M}(|\nabla\times u|^{q_0})>N_1^{q_0}\}\cap K_r(\hat{x},\hat{t})|\geq\varepsilon|K_r(\hat{x},\hat{t})|,
\end{eqnarray*}
then
\begin{eqnarray*}
\begin{aligned}
\lefteqn{K_r(\hat{x},\hat{t})\cap Q_1\subset\{(x,t)\in Q_1:\mathcal{M}(|\nabla\times u|^{q_0})>1\}\cup\{(x,t)\in Q_1:\Theta>1\}}\hspace*{28mm}\\
&\ \cup\{(x,t)\in Q_1:\mathcal{M}(|f|^{q_0})>\tilde{\delta}^{q_0}\}\cup\{(x,t)\in Q_1:\mathcal{M}(|g|^{q_0})>1\},
\end{aligned}
\end{eqnarray*}
Proposition 3.2 follows from Lemma 2.4 directly.
\qed

\bigskip

Applying Proposition 3.2 to
$\frac{u}{N^{k-1}_1}$,~~$\frac{p}{N^{k-1}_1}$,~~$\frac{f}{N^{k-1}_1}$~~ \mbox{and} ~~$\frac{g}{N^{k-1}_1}$, we immediately have:

\bigskip

\begin{corollary}
Under the assumptions of Proposition 3.2, we have
\begin{eqnarray}\label{d4.41}
\begin{aligned}
\lefteqn{|\{(x,t)\in Q_1:\mathcal{M}(|\nabla\times u|^{q_0})>N^{q_0k}_1\}|}\hspace*{5mm}\\
\leq&\ 10^{n+2}\varepsilon\Big(|\{(x,t)\in Q_1:\mathcal{M}(|\nabla\times u|^{q_0})>N^{q_0(k-1)}_1\}|\\
&\mbox{}\hskip0.6cm +
|\{(x,t)\in Q_1:\Theta>N^{q_0(k-1)}_1\}|+|\{(x,t)\in Q_1:\mathcal{M}(|f|^{q_0})>\tilde{\delta}^{q_0}N^{q_0(k-1)}_1\}|\\
&\mbox{}\hskip1.6cm +|\{(x,t)\in Q_1:\mathcal{M}(|g|^{q_0})>N^{q_0(k-1)}_1\}|\Big)
\end{aligned}
\end{eqnarray}
for any $k\in\mathbb{N}_+$.
\end{corollary}

\bigskip

To do our iteration for improving the regularity of $Du$,
we set
\begin{eqnarray}\label{d4.50}
p_{j+1}=\frac{pp_{j}}{p(1-\gamma)+p_{j}\gamma}~~\mbox{for}~~j=0,1,2,\ldots,
\end{eqnarray}
where $1<p_0<p$ are given in Theorem 3.1 and $\gamma>0$ is the constant in hypothesis ($H_2'$). Observe that $\{p_{j}\}$ is increasing and $\lim\limits_{j\rightarrow\infty}p_j=p$.

\bigskip

We first give an estimate of $\Theta$.
\begin{lemma}
Let $1<p_0<p<\infty$, $q_0=\frac{1+p_0}{2}$, and $\Theta$ and $\{p_{j}\}$ be given by (\ref{d4.37}) and (\ref{d4.50}) respectively. Suppose $u\in L_{p}(Q_{7};\mathbb{R}^{n})$ and $Du\in L_{p_j}(Q_{7};\mathbb{R}^{n^2})$. Then
we have $\Theta\in L_{p_{j+1}/q_0}(Q_{1})$ and
\begin{eqnarray}\label{d8.50}
\|\Theta\|^{p_{j+1}/q_0}_{L_{p_{j+1}/q_0}(Q_{1})}
\leq C\Big(\|\mathcal{M}(|Du|^{q_0})\|^{p_{j}/q_0}_{L_{p_{j}/q_0}(Q_{1})}
+\|\mathcal{M}(|u|^{q_0})\|^{p/q_0}_{L_{p/q_0}(Q_{1})}\Big),
\end{eqnarray}
where $C$ depends only on $n$, $p_0$ and $p$.
\end{lemma}
\vskip 2mm\noindent{\bf Proof.}\quad
By (\ref{d4.37}) and Poincar\'{e}'s inequality, we have
\begin{eqnarray*}
\begin{aligned}
\Theta
=&\  \sup_{r>0}\Bigg\{\Big(r^{-q_0}\fint_{K_r(x,t)\cap Q_7}|u-\bar{u}_{B_r(x)}|^{q_0}dyds\Big)^{1-\gamma}
 \Big(\fint_{K_r(x,t)\cap Q_7}|u-\bar{u}_{B_r(x)}|^{q_0}dyds\Big)^{\gamma}\Bigg\}\\
\leq&\ C\sup_{r>0}\Bigg\{\Big(\fint_{K_r(x,t)\cap Q_7}|Du|^{q_0}dyds\Big)^{1-\gamma}
\Big(\fint_{K_r(x,t)\cap Q_7}|u|^{q_0}dyds\Big)^{\gamma}\Bigg\}\\
\leq&\ C\big(\mathcal{M}(|Du|^{q_0})\big)^{1-\gamma}
\big(\mathcal{M}(|u|^{q_0})\big)^{\gamma},
\end{aligned}
\end{eqnarray*}
where $C$ depends only on $n$ and $p_0$.
From Young's inequality and (\ref{d4.50}), it follows that
\begin{eqnarray*}
\begin{aligned}
\int_{Q_1}|\Theta|^{\frac{p_{j+1}}{q_0}}dxdt
\leq&\ C\int_{Q_1}\Bigg\{\big(\mathcal{M}(|Du|^{q_0})\big)^{\frac{p_{j+1}(1-\gamma)}{q_0}}
\big(\mathcal{M}(|u|^{q_0})\big)^{\frac{p_{j+1}\gamma}{q_0}}\Bigg\}dxdt\\
\leq&\ C\Bigg\{\int_{Q_1}\big(\mathcal{M}(|Du|^{q_0})\big)^{\frac{p_{j}}{q_0}}dxdt +\int_{Q_1}\big(\mathcal{M}(|u|^{q_0})\big)^{\frac{p}{q_0}}dxdt\Bigg\},
\end{aligned}
\end{eqnarray*}
where $C$ depends only on $n$, $p_0$ and $p$.
Thus, (\ref{d8.50}) holds.
\qed

\bigskip

\begin{proposition}
Let $0<\varepsilon<1<p_0<p<\infty$, $q_0=\frac{1+p_0}{2}$, $\{p_{j}\}$ be given by (\ref{d4.50}) and $u$, $\pi$, $f$ and $g$ be given by Theorem 3.1. Assume that $Q_7\subset \Omega_{T}$, $(H_1)$ and $(H_2')$ hold with some $0<\delta<1$ depending only on $\varepsilon$, $n$, $\lambda$, $p_0$ and $p$.
If $Du\in L_{p_j}(Q_{7};\mathbb{R}^{n^2})$, then $Du\in L_{p_{j+1}}(Q_{\frac{1}{2}};\mathbb{R}^{n^2})$ and
\begin{eqnarray}\label{d10.49}
\begin{aligned}
\|Du\|^{p_{j+1}}_{L_{p_{j+1}}(Q_{\frac{1}{2}})}
\leq C\varepsilon\Big(\|Du\|^{p_{j}}_{L_{p_{j}}(Q_{7})}
+\|Du\|^{p_{j+1}}_{L_{p_{j}}(Q_{7})}\Big)+\tilde{C}M_j,
\end{aligned}
\end{eqnarray}
where $C$ depends only on $n$, $\lambda$, $p_0$ and $p$, $\tilde{C}$ depends only on $\varepsilon$, $n$, $\lambda$, $p_0$ and $p$ and
\begin{eqnarray*}
M_j=\|u\|^{p}_{L_{p}(Q_{7})}
+\|u\|^{p_{j+1}}_{L_{p}(Q_{7})}+\|f\|^{p_{j+1}}_{L_{p}(Q_{7})}
+\|g\|^{p_{j+1}}_{L_{p}(Q_{7})}.
\end{eqnarray*}
\end{proposition}
\vskip 2mm\noindent{\bf Proof.}\quad
Let $N_1$ and $\tilde{\delta}$ be given by Proposition 3.2. Take $\delta=\varepsilon\tilde{\delta}$.
Denote
\begin{eqnarray}\label{n.33}
\tilde{u}=\frac{u}{N}, ~~\tilde{\pi}=\frac{\pi}{N},~~\tilde{f}=\frac{ f}{N}~~\mbox{and}~~\tilde{g}=\frac{g}{N},
\end{eqnarray}
where
\begin{eqnarray*}
N=\bar{C}^{\frac{1}{q_0}}N_1^{-1}\varepsilon^{-\frac{1}{q_0}}\Big(\|u\|_{L_p(Q_{7})}+\|f\|_{L_p(Q_{7})}
+\|g\|_{L_p(Q_{7})}+\delta\|Du\|_{L_{q_0}(Q_{7})}\Big)
\end{eqnarray*}
with the constant $\bar{C}>0$ given by the following (\ref{n1.74}) and depending only on $n$, $\lambda$, $p_0$ and $p$.
It is easy to see that $\tilde{u}$, $\tilde{\pi}$, $\tilde{f}$ and $\tilde{g}$ also satisfy (\ref{s0.0}). Then, by Lemma 2.3 (ii), we have
\begin{eqnarray}\label{n.74}
\begin{aligned}
|\{(x,t)\in Q_1:\mathcal{M}(|\nabla\times \tilde{u}|^{q_0})>N^{q_0}_1\}|
\leq \frac{C}{N^{q_0}_1}\|\nabla\times \tilde{u}\|^{q_0}_{L_{q_0}(Q_{2})}.
\end{aligned}
\end{eqnarray}

Similar to (\ref{n.78}), we obtain
\begin{equation*}
(\nabla\times \tilde{u})_t-\mbox{div}\big(\overline{A}_{B_7}(t)\nabla(\nabla\times \tilde{u})\big)
=\nabla\times\mbox{div}~\Big(\big(A-\overline{A}_{B_7}(t)\big)D\tilde{u}+\tilde{f}\Big)~~\mbox{in}~~Q_{7}.
\end{equation*}
By virtue of (\ref{r.5}) and $(H'_2)$,
\begin{equation*}
\begin{aligned}
\lefteqn{\|\nabla\times \tilde{u}\|_{L_{q_0}(Q_{2})}}\hspace*{6mm}\\
\leq&\ C\Bigg\{\|\big(A-\overline{A}_{B_{7}}(t)\big)D\tilde{u}\|_{L_{q_0}(Q_{7})}+\|\tilde{f}\|_{L_{q_0}(Q_{7})} +\Big(\int_{-49}^0\|\nabla\times \tilde{u}(\cdot,t)\|^{q_0}_{W^{-1,q_0}(B_{7})}dt\Big)^{\frac{1}{q_0}}\Bigg\}\\
\leq&\ C\Big(\delta\|D\tilde{u}\|_{L_{q_0}(Q_{7})}+\|\tilde{f}\|_{L_{q_0}(Q_{7})}+\|\tilde{u}\|_{L_{q_0}(Q_{7})}\Big),
\end{aligned}
\end{equation*}
where $C$ depends only on $n$, $\lambda$ and $p_0$.
This and (\ref{n.74}) imply
\begin{eqnarray}\label{n1.74}
\begin{aligned}
\lefteqn{|\{(x,t)\in Q_1:\mathcal{M}(|\nabla\times \tilde{u}|^{q_0})>N^{q_0}_1\}|}\hspace*{26mm}\\
\leq&\ \frac{C}{N^{q_0}_1}\Big(\delta\|D\tilde{u}\|_{L_{q_0}(Q_{7})}+\|\tilde{f}\|_{L_{q_0}(Q_{7})}+\|\tilde{u}\|_{L_{q_0}(Q_{7})}\Big)^{q_0}\\
\leq&\ \frac{\bar{C}}{(N_1N)^{q_0}}\Big(\delta\|Du\|_{L_{q_0}(Q_{7})}+\|f\|_{L_{p}(Q_{7})}
+\|u\|_{L_{p}(Q_{7})}\Big)^{q_0}
\leq \varepsilon|Q_1|.
\end{aligned}
\end{eqnarray}

Hence, Corollary 3.5 holds for $\tilde{u}$, $\tilde{\pi}$, $\tilde{f}$ and $\tilde{g}$. Multiplying (\ref{d4.41}) by $N^{kp_{j+1}}_1$ and taking sum in $k$, one has
\begin{equation*}
\begin{aligned}
\lefteqn{\sum^\infty_{k=1}N^{kp_{j+1}}_1|\{(x,t)\in Q_1:\mathcal{M}(|\nabla\times \tilde{u}|^{q_0})>N^{q_0k}_1\}|}\hspace*{8mm}\\
\leq&\ 10^{n+2}\varepsilon N^{p_{j+1}}_1\sum^\infty_{k=1}\Big(N^{(k-1)p_{j+1}}_1|\{(x,t)\in Q_1:\mathcal{M}(|\nabla\times \tilde{u}|^{q_0})>N^{q_0(k-1)}_1\}|\\
&\mbox{}\hskip2cm +N^{(k-1)p_{j+1}}_1|\{(x,t)\in Q_1:\tilde{\Theta}>N^{q_0(k-1)}_1\}|\\
&\mbox{}\hskip3cm+N^{(k-1)p_{j+1}}_1|\{(x,t)\in Q_1:\frac{1}{\tilde{\delta}^{q_0}}\mathcal{M}(|\tilde{f}|^{q_0})>N^{q_0(k-1)}_1\}|\\
&\mbox{}\hskip4cm+N^{(k-1)p_{j+1}}_1|\{(x,t)\in Q_1:\mathcal{M}(|\tilde{g}|^{q_0})>N^{q_0(k-1)}_1\}|\Big),
\end{aligned}
\end{equation*}
where $\tilde{\Theta}$ is defined by (\ref{d4.37}) with $u$ replaced by $\tilde{u}$.
From Lemma 3.6, we obtain $\tilde{\Theta}\in L_{p_{j+1}/q_0}(Q_{1})$. Then, it follows from Lemma 2.1, Lemma 2.3 (i) and (\ref{n.33}) that
\begin{eqnarray}\label{n.34}
\begin{aligned}
\lefteqn{\|\mathcal{M}(|\nabla\times \tilde{u}|^{q_0})\|^{p_{j+1}/q_0}_{L_{p_{j+1}/q_0}(Q_{1})}}\hspace*{18mm}\\
\leq&\ C\varepsilon\Big(\|\mathcal{M}(|\nabla\times \tilde{u}|^{q_0})\|^{p_{j+1}/q_0}_{L_{p_{j+1}/q_0}(Q_{1})}
+\|\tilde{\Theta}\|^{p_{j+1}/q_0}_{L_{p_{j+1}/q_0}(Q_{1})}\Big)+\frac{C\varepsilon^{\frac{p_{j+1}}{q_0}}}{\tilde{\delta}^{p_{j+1}}},
\end{aligned}
\end{eqnarray}
where we use the fact
\begin{eqnarray*}
\begin{aligned}
\|\frac{1}{\tilde{\delta}^{q_0}}\mathcal{M}(|\tilde{f}|^{q_0})\|^{p_{j+1}/q_0}_{L_{p_{j+1}/q_0}(Q_{1})}
\leq \frac{C}{\tilde{\delta}^{p_{j+1}}}\|\tilde{f}\|^{p_{j+1}}_{L_{p_{j+1}}(Q_{7})}
\leq \frac{C\varepsilon^{\frac{p_{j+1}}{q_0}}}{\tilde{\delta}^{p_{j+1}}},
\end{aligned}
\end{eqnarray*}
\begin{eqnarray*}
\begin{aligned}
\|\mathcal{M}(|\tilde{g}|^{q_0})\|^{p_{j+1}/q_0}_{L_{p_{j+1}/q_0}(Q_{1})}
\leq C\|\tilde{g}\|^{p_{j+1}}_{L_{p_{j+1}}(Q_{7})}
\leq C\varepsilon^{\frac{p_{j+1}}{q_0}},
\end{aligned}
\end{eqnarray*}
and $C$ depends only on $n$, $\lambda$, $p_0$ and $p$. Take $\varepsilon\leq\frac{1}{2C}$ with the same $C$ in (\ref{n.34}) and then
\begin{eqnarray*}
\begin{aligned}
\|\mathcal{M}(|\nabla\times \tilde{u}|^{q_0})\|^{p_{j+1}/q_0}_{L_{p_{j+1}/q_0}(Q_{1})}
\leq&\ C\varepsilon\|\tilde{\Theta}\|^{p_{j+1}/q_0}_{L_{p_{j+1}/q_0}(Q_{1})}
+\frac{C\varepsilon^{\frac{p_{j+1}}{q_0}}}{\tilde{\delta}^{p_{j+1}}},
\end{aligned}
\end{eqnarray*}
where $C$ depends only on $n$, $\lambda$, $p_0$ and $p$.

Rewriting the above inequality in terms of $u$, we have
\begin{eqnarray*}
\begin{aligned}
\|\nabla\times u\|^{p_{j+1}}_{L_{p_{j+1}}(Q_{1})}
\leq&\ C\varepsilon\|\Theta\|^{p_{j+1}/q_0}_{L_{p_{j+1}/q_0}(Q_{1})}
+\frac{C}{\tilde{\delta}^{p_{j+1}}}\Big(\|u\|^{p_{j+1}}_{L_{p}(Q_{7})}+\|f\|^{p_{j+1}}_{L_{p}(Q_{7})}\\
&\mbox{}\hskip5cm
+\|g\|^{p_{j+1}}_{L_{p}(Q_{7})}+\delta^{p_{j+1}}\|Du\|^{p_{j+1}}_{L_{q_0}(Q_{7})}\Big).
\end{aligned}
\end{eqnarray*}
Since $\delta=\varepsilon\tilde{\delta}$, then
\begin{eqnarray*}
\begin{aligned}
\|\nabla\times u\|^{p_{j+1}}_{L_{p_{j+1}}(Q_{1})}
\leq&\ C\varepsilon\Big(\|\Theta\|^{p_{j+1}/q_0}_{L_{p_{j+1}/q_0}(Q_{1})}+\|Du\|^{p_{j+1}}_{L_{p_j}(Q_{7})}\Big)\\
&\ +\tilde{C}\Big(\|u\|^{p_{j+1}}_{L_{p}(Q_{7})}+\|f\|^{p_{j+1}}_{L_{p}(Q_{7})}
+\|g\|^{p_{j+1}}_{L_{p}(Q_{7})}\Big),
\end{aligned}
\end{eqnarray*}
where $C$ depends only on $n$, $\lambda$, $p_0$ and $p$ and $\tilde{C}$ depends only on $\varepsilon$, $n$, $\lambda$, $p_0$ and $p$. (Observe that $\tilde{\delta}$ depends on $\varepsilon$.)
By (\ref{d8.50}) and Lemma 2.3 (i), it leads to
\begin{eqnarray}\label{n.35}
\begin{aligned}
\|\nabla\times u\|^{p_{j+1}}_{L_{p_{j+1}}(Q_{1})}
\leq&\ C\varepsilon\Big(\|Du\|^{p_{j}}_{L_{p_{j}}(Q_{7})}+\|Du\|^{p_{j+1}}_{L_{p_j}(Q_{7})}\Big)\\
&\ +\tilde{C}\Big(\|u\|^{p}_{L_{p}(Q_{7})}+\|u\|^{p_{j+1}}_{L_{p}(Q_{7})}+\|f\|^{p_{j+1}}_{L_{p}(Q_{7})}
+\|g\|^{p_{j+1}}_{L_{p}(Q_{7})}\Big),
\end{aligned}
\end{eqnarray}
where $C$ depends only on $n$, $\lambda$, $p_0$ and $p$ and $\tilde{C}$ depends only on $\varepsilon$, $n$, $\lambda$, $p_0$ and $p$.

From $\mbox{div}~u=g$, it follows from (\ref{n.77}) that
\begin{eqnarray*}
\begin{aligned}
\Delta u
=\mbox{div}~\nabla\times u+\nabla g.
\end{aligned}
\end{eqnarray*}
Then, using $W^{1,p}$ estimates to $u$,
we obtain
\begin{eqnarray*}
\begin{aligned}
\|Du\|^{p_{j+1}}_{L_{p_{j+1}}(Q_{\frac{1}{2}})}
\leq&\ C\Big(\|\nabla\times u\|^{p_{j+1}}_{L_{p_{j+1}}(Q_{1})}
+\|u\|^{p_{j+1}}_{L_{p_{j+1}}(Q_{1})}
+\|g\|^{p_{j+1}}_{L_{p_{j+1}}(Q_{1})}\Big),
\end{aligned}
\end{eqnarray*}
where $C$ depends only on $n$, $p_0$ and $p$.
Substituting (\ref{n.35}) into it, we get
\begin{eqnarray*}
\begin{aligned}
\|Du\|^{p_{j+1}}_{L_{p_{j+1}}(Q_{\frac{1}{2}})}
\leq&\ C\varepsilon\Big(\|Du\|^{p_{j}}_{L_{p_{j}}(Q_{7})}+\|Du\|^{p_{j+1}}_{L_{p_j}(Q_{7})}\Big)\\
&\ +\tilde{C}\Big(\|u\|^{p}_{L_{p}(Q_{7})}
+\|u\|^{p_{j+1}}_{L_{p}(Q_{7})}+\|f\|^{p_{j+1}}_{L_{p}(Q_{7})}
+\|g\|^{p_{j+1}}_{L_{p}(Q_{7})}\Big),
\end{aligned}
\end{eqnarray*}
where $C$ depends only on $n$, $\lambda$, $p_0$ and $p$ and $\tilde{C}$ depends only on $\varepsilon$, $n$, $\lambda$, $p_0$ and $p$.
\qed

\bigskip

\begin{corollary}
Under the assumptions of Proposition 3.7,
if $Du\in L_{p_j}(Q_{r};\mathbb{R}^{n^2})$, then $Du\in L_{p_{j+1}}(Q_{\tau r};\mathbb{R}^{n^2})$ and
\begin{eqnarray}\label{n.138}
\begin{aligned}
\|Du\|^{p_{j+1}}_{L_{p_{j+1}}(Q_{\tau r})}
\leq\frac{C\varepsilon}{(r-\tau r)^{(n+2)(\frac{p_{j+1}}{p_j}-1)}}\|Du\|^{p_{j+1}}_{L_{p_{j}}(Q_{r})}
+\frac{\tilde{C}}{(r-\tau r)^p}M,
\end{aligned}
\end{eqnarray}
where
$0<\tau<1$, $0< r\leq7$, $C$ depends only on $n$, $\lambda$, $p_0$ and $p$, $\tilde{C}$ depends only on $\varepsilon$, $n$, $\lambda$, $p_0$ and $p$, and
\begin{eqnarray}\label{n.38}
M=\|u\|^{p}_{L_{p}(Q_{7})}+\|f\|^{p}_{L_{p}(Q_{7})}
+\|g\|^{p}_{L_{p}(Q_{7})}+1.
\end{eqnarray}
\end{corollary}
\vskip 2mm\noindent{\bf Proof.}\quad
By virtue of Proposition 3.7 and standard scaling analysis, we get
\begin{eqnarray*}
\begin{aligned}
\|Du\|^{p_{j+1}}_{L_{p_{j+1}}(Q_{\frac{r}{2}})}
\leq&\ C\varepsilon\Big(\|Du\|^{p_{j}}_{L_{p_{j}}(Q_{7r})}
+\frac{1}{r^{(n+2)(\frac{p_{j+1}}{p_j}-1)}}\|Du\|^{p_{j+1}}_{L_{p_{j}}(Q_{7r})}\Big)+\frac{\tilde{C}}{r^p}M_j\\
\leq&\
\frac{C\varepsilon}{r^{(n+2)(\frac{p_{j+1}}{p_j}-1)}}\|Du\|^{p_{j+1}}_{L_{p_{j}}(Q_{7r})}+\frac{\tilde{C}}{r^p}M
\end{aligned}
\end{eqnarray*}
for $0< r\leq1$, where
$\|Du\|^{p_{j}}_{L_{p_{j}}(Q_{7r})}\leq\|Du\|^{p_{j+1}}_{L_{p_{j}}(Q_{7r})}+1~~\mbox{and}~~M_j+1\leq 4M$ are used. Applying standard scaling and covering arguments again, (\ref{n.138}) is obtained.
\qed

\bigskip

We are now ready to prove Theorem 3.1.

\bigskip

\noindent{\bf Proof of Theorem 3.1.}\quad
For any integer $j_0\geq1$, set
\begin{eqnarray*}
r_j=1-\big(\frac{1}{2}\big)^{j_0-j+1},~~j=0,1,2,...,j_0.
\end{eqnarray*}
It is clear that  $r_{j_0}=\frac{1}{2}$ and $\frac{1}{2}\leq r_{j+1}<r_j<1$.
In Corollary 3.8, set $r=r_j$, $\tau=\frac{r_{j+1}}{r_j}$ and then
\begin{eqnarray*}
r-\tau r=r_j-r_{j+1}=2^{-(j_0-j+1)}.
\end{eqnarray*}
By (\ref{d4.50}), we get
$\frac{p_{j+1}}{p_j}-1\leq \frac{p}{p_0}(1-\gamma)^j$ and then
\begin{eqnarray*}
(r-\tau r)^{-(n+2)(\frac{p_{j+1}}{p_j}-1)}=(r_j-r_{j+1})^{-(n+2)(\frac{p_{j+1}}{p_j}-1)}
\leq 2^{C_0(j_0-j+1)(1-\gamma)^j}
\end{eqnarray*}
with $C_0=(n+2)\frac{p}{p_0}$.
Therefore, if $Du\in L_{p_j}(Q_{r_j};\mathbb{R}^{n^2})$, by Corollary 3.8, we obtain $Du\in L_{p_{j+1}}(Q_{r_{j+1}};\mathbb{R}^{n^2})$ and
\begin{eqnarray}\label{n5.35}
\begin{aligned}
\|Du\|^{p_{j+1}}_{L_{p_{j+1}}(Q_{r_{j+1}})}
\leq C\varepsilon2^{C_0(j_0-j+1)(1-\gamma)^j}\|Du\|^{p_{j+1}}_{L_{p_{j}}(Q_{r_j})}
+2^{(j_0-j+1)p}\tilde{C}M
\end{aligned}
\end{eqnarray}
for $j=0,1,2,...,j_0-1$, where $C$ depends only on $n$, $\lambda$, $p_0$ and $p$, $\tilde{C}$ depends only on $\varepsilon$, $n$, $\lambda$, $p_0$ and $p$, and $M$ is given by (\ref{n.38}).

Using (\ref{n5.35}) and H\"{o}lder's inequality\footnote{$(a+b)^p\leq2^{p-1}(a^p+b^p)$~for~any~$a,b>0$~and~$p>1$.} iteratively from $j=j_0-1$ to $j=0$, one has
\begin{eqnarray*}
\begin{aligned}
\|Du\|^{p_{j_0}}_{L_{p_{j_0}}(Q_{r_{j_0}})}
\leq&\ 2^{\alpha_{0}+C_0\beta_{0}}
(C\varepsilon)^{\zeta_{0}+1}
\|Du\|^{p_{j_0}}_{L_{p_{0}}(Q_{r_0})}\\
&\ +
\sum_{j=1}^{j_0}2^{\alpha_{j}+C_0\beta_{j}+\frac{(j_0-j+2)pp_{j_0}}{p_{j}}}
(C\varepsilon)^{\zeta_{j}}(\tilde{C}M)^\frac{p_{j_0}}{p_{j}},
\end{aligned}
\end{eqnarray*}
where
\begin{eqnarray*}
\zeta_{0}=\sum_{j=1}^{j_0-1}\frac{p_{j_0}}{p_{j}},
~~\zeta_{j}=\sum_{k=j+1}^{j_0}\frac{p_{j_0}}{p_{k}},
\end{eqnarray*}
\begin{eqnarray*}
\alpha_{0}=\sum_{j=1}^{j_0-1}\Big(\frac{p_{j_0}}{p_{j}}-1\Big),
~~\alpha_{j}=\sum_{k=j}^{j_0-1}\Big(\frac{p_{j_0}}{p_{k}}-1\Big),
\end{eqnarray*}
\begin{eqnarray*}
\beta_{0}=\sum_{j=1}^{j_0}(j_0-j+2)(1-\gamma)^{j-1}\frac{p_{j_0}}{p_{j}}
~~\mbox{and}~~\beta_{j}=\sum_{k=j}^{j_0-1}(j_0-k+1)(1-\gamma)^{k}\frac{p_{j_0}}{p_{k+1}}
\end{eqnarray*}
for $j=1,2,...,j_0$. Observe $Du\in L_{p_0}(Q_7;\mathbb{R}^{n^2})$. Since
\begin{eqnarray*}
\zeta_{0}\geq j_0-1,~~\zeta_{j}\geq j_0-j,~~\alpha_{0}\leq\frac{p}{p_0}j_0,~~\alpha_{j}\leq\frac{p}{p_0}(j_0-j),
\end{eqnarray*}
\begin{eqnarray*}
\begin{aligned}
\beta_{0}\leq 2j_0\frac{p}{p_0}\sum_{j=1}^{j_0}(1-\gamma)^{j-1}
\leq \frac{2p}{\gamma p_0}j_0,
\end{aligned}
\end{eqnarray*}
\begin{eqnarray*}
\begin{aligned}
\beta_{j}\leq&\ 2(j_0-j)\frac{p}{p_0}\sum_{k=j}^{j_0-1}(1-\gamma)^{k}
\leq \frac{2p}{\gamma p_0}(j_0-j),
\end{aligned}
\end{eqnarray*}
\begin{eqnarray*}
(\tilde{C}M)^\frac{p_{j_0}}{p_{j}}
\leq(\tilde{C}M)^\frac{p}{p_{0}}
~~\mbox{and}~~
\frac{(j_0-j+2)pp_{j_0}}{p_{j}}\leq\frac{(j_0-j+2)p^2}{p_0}
\end{eqnarray*}
for $j=1,2,...,j_0$, we take $\varepsilon>0$ sufficiently small so that $2^{\frac{2p^2}{p_0}+\frac{2C_0p}{\gamma p_0}}C\varepsilon<1$ and then
\begin{eqnarray*}
\begin{aligned}
\|Du\|^{p_{j_0}}_{L_{p_{j_0}}(Q_{r_{j_0}})}
\leq&\ \big(2^{\frac{p}{p_{0}}+\frac{2C_0p}{\gamma p_0}}C\varepsilon\big)^{j_0}\|Du\|^{p_{j_0}}_{L_{p_{0}}(Q_{1})}
+4^{\frac{p^2}{p_0}}(\tilde{C}M)^\frac{p}{p_{0}}
\sum_{j=1}^{j_0}\big(2^{\frac{2p^2}{p_0}+\frac{2C_0p}{\gamma p_0}}C\varepsilon\big)^{j_0-j}\\
\leq&\ \big(2^{\frac{p}{p_{0}}+\frac{2C_0p}{\gamma^2 p_0}}C\varepsilon\big)^{j_0}\|Du\|^{p_{j_0}}_{L_{p_{0}}(Q_{1})}+
\frac{4^{\frac{p^2}{p_0}}(\tilde{C}M)^\frac{p}{p_{0}}}{1-2^{\frac{2p^2}{p_0}+\frac{2C_0p}{\gamma p_0}}C\varepsilon}.
\end{aligned}
\end{eqnarray*}
Observe that as $\varepsilon$ is chosen, $\tilde{C}$ is determined and depends only on $n$, $\lambda$, $\gamma$, $p_0$ and $p$.

Let $j_0\rightarrow\infty$ and since $\lim\limits_{j_0\rightarrow\infty}p_{j_0}=p$ and $2^{\frac{p}{p_{0}}+\frac{2C_0p}{\gamma^2 p_0}}C\varepsilon<1$, we deduce
\begin{eqnarray*}
\begin{aligned}
\|Du\|^{p}_{L_{p}(Q_{\frac{1}{2}})}
\leq CM^{\frac{p}{p_0}},
\end{aligned}
\end{eqnarray*}
where $C$ depends only on $n$, $\lambda$, $\gamma$, $p_0$ and $p$.
This and (\ref{n.38}) imply
\begin{eqnarray}\label{n.39}
\begin{aligned}
\|Du\|_{L_{p}(Q_{\frac{1}{2}})}
\leq&\ C\Big(\big(\|u\|_{L_{p}(Q_{7})}+\|f\|_{L_{p}(Q_{7})}
+\|g\|_{L_{p}(Q_{7})}\big)^{\frac{p}{p_0}}+1\Big),
\end{aligned}
\end{eqnarray}
where $C$ depends only on $n$, $\lambda$, $\gamma$, $p_0$ and $p$.

Finally, we use the scaling arguments to modify (\ref{n.39}). Indeed,
for any constant $H>0$,
using (\ref{n.39}) to $uH$, $\pi H$, $fH$ and $gH$ and dividing both sides by $H$, we have
\begin{eqnarray*}
\begin{aligned}
\|Du\|_{L_{p}(Q_{\frac{1}{2}})}
\leq C\Big(H^{\frac{p}{p_0}-1}\big(\|u\|_{L_{p}(Q_{7})}+\|f\|_{L_{p}(Q_{7})}
+\|g\|_{L_{p}(Q_{7})}\big)^{\frac{p}{p_0}}+\frac{1}{H}\Big).
\end{aligned}
\end{eqnarray*}
Taking
$H=\big(\|u\|_{L_{p}(Q_{7})}+\|f\|_{L_{p}(Q_{7})}
+\|g\|_{L_{p}(Q_{7})}\big)^{-1},
$
the above inequality becomes
\begin{eqnarray*}
\begin{aligned}
\|Du\|_{L_{p}(Q_{\frac{1}{2}})}
\leq C\Big(\|u\|_{L_{p}(Q_{7})}+\|f\|_{L_{p}(Q_{7})}
+\|g\|_{L_{p}(Q_{7})}\Big),
\end{aligned}
\end{eqnarray*}
where $C$ depends only on $n$, $\lambda$, $\gamma$, $p_0$ and $p$.
That is, (\ref{n1}) is obtained.
\qed

\bigskip

\section{Regularity for Stokes systems in non-divergence form}\setcounter{equation}{0}
This section is devoted to proving Theorem 1.2, the interior $L_{p}$ regularity for Stokes systems in non-divergence form. Recall $(H_2')$ is given in section 3. By standard scaling and covering arguments, we only need to show:
\begin{theorem}
Let $1<p_0<p<\infty$. Suppose $u\in W_{p_0}^{2,1}(\Omega_{T};\mathbb{R}^n)\cap L_p(\Omega_{T};\mathbb{R}^n)$, $\nabla \pi\in L_{1}(\Omega_{T};\mathbb{R}^n)$, $f\in L_{p}(\Omega_{T};\mathbb{R}^{n})$ and $Dg\in L_{p}(\Omega_{T};\mathbb{R}^{n})$ satisfying (\ref{s0.1}). Assume that $Q_7\subset\Omega_T$, $(H_1)$ and $(H_2')$ hold with some $0<\delta<1$ depending only on $n$, $\lambda$, $\gamma$, $p_0$ and $p$.
Then we have $D^2u\in L_{p}(Q_\frac{1}{2};\mathbb{R}^{n^3})$ and
\begin{eqnarray}\label{n2}
\|D^2u\|_{L_p(Q_\frac{1}{2})}
\leq
C\Big(\|u\|_{L_p(Q_7)}+\|f\|_{L_p(Q_7)}+\|Dg\|_{L_p(Q_7)}\Big),
\end{eqnarray}
where $C$ depends only on $n$, $\lambda$, $\gamma$, $p_0$ and $p$.
\end{theorem}

\bigskip

The proof of Theorem 4.1 is very similar to that of Theorem 3.1 and we only give a sketch of it. To establish the iteration inequality (see (\ref{d6.48})), the main step is also to obtain the proposition below, which is a consequence of Lemma 4.4. Recall that we abbreviate $\mathcal{M}_{Q_7}(f)=\mathcal{M}(f)$ for any function $f$.

\bigskip

\begin{proposition}
Set $0<\varepsilon<1<p_0<\infty$ and $q_0=\frac{1+p_0}{2}$. Suppose $u\in W_{p_0}^{2,1}(\Omega_{T};\mathbb{R}^{n})$, $\nabla \pi\in L_{1}(\Omega_{T};\mathbb{R}^{n})$, $f\in L_{p_0}(\Omega_{T};\mathbb{R}^{n})$, $Dg\in L_{p_0}(\Omega_{T};\mathbb{R}^{n})$ satisfying (\ref{s0.1}), $Q_7\subset\Omega_{T}$ and $(H_1)$ holds. There exists a constant $0<\tilde{\delta}<1$ depending only on $\varepsilon$, $n$, $\lambda$ and $p_0$. Assume that $(H_2')$ holds with $0<\delta\leq\tilde{\delta}$.  For any $(x,t)\in Q_1$, let
\begin{eqnarray}\label{d2.37}
\Theta(x,t):=\sup_{r>0}r^{q_0(\gamma-1)}\fint_{K_r(x,t)\cap Q_7}|Du(y,s)-\overline{Du}_{B_r(x)}(s)|^{q_0}dyds,
\end{eqnarray}
where $\gamma$ is given in ($H_2'$).
There exists a constant $N_1>1$ depending only on $n$, $\lambda$ and $p_0$ such that if
\begin{eqnarray*}
|\{(x,t)\in Q_1:\mathcal{M}(|D\nabla\times u|^{q_0})>N_1^{q_0}\}|\leq \varepsilon|Q_1|,
\end{eqnarray*}
then
\begin{equation*}
\begin{aligned}
\lefteqn{|\{(x,t)\in Q_1:\mathcal{M}(|D\nabla\times u|^{q_0})>N_1^{q_0}\}|}\hspace*{16mm}\\
\leq&\ 10^{n+2}\varepsilon\Big(|\{(x,t)\in Q_1:\mathcal{M}(|D\nabla\times u|^{q_0})>1\}|+|\{(x,t)\in Q_1:\Theta>1\}|\\
& +|\{(x,t)\in Q_1:\mathcal{M}(|f|^{q_0})>\tilde{\delta}^{q_0}\}|+|\{(x,t)\in Q_1:\mathcal{M}(|Dg|^{q_0})>1\}|\Big).
\end{aligned}
\end{equation*}
\end{proposition}

\bigskip

\begin{remark}
The scaling of $\Theta$ is between $Du$ and $D^2u$. In the proof of Proposition 4.2, $\Theta$ will be applied to control $Du$. In Lemma 4.6, $\Theta$ will be estimated by $u$ and $D^2u$.
\end{remark}

\bigskip

\begin{lemma}
Under the assumptions of Proposition 4.2, there exists a constant $N_1>1$ depending only on $n$, $\lambda$ and $p_0$ such that if
\begin{eqnarray}\label{d2.20}
\begin{aligned}
\lefteqn{\{(x,t)\in Q_1:\mathcal{M}(|D\nabla\times u|^{q_0})\leq1\}\cap\{(x,t)\in Q_1:\Theta\leq1\}}\hspace*{2mm}\\
&\cap\{(x,t)\in Q_1:\mathcal{M}(|f|^{q_0})\leq\tilde{\delta}^{q_0}\}
\cap\{(x,t)\in Q_1:\mathcal{M}(|Dg|^{q_0})\leq1\}\cap K_r(\hat{x},\hat{t})
\neq\emptyset
\end{aligned}
\end{eqnarray}
for some $(\hat{x},\hat{t})\in Q_1$ and $0<r\leq1$, then
\begin{eqnarray}\label{d9.20}
|\{(x,t)\in Q_1:\mathcal{M}(|D\nabla\times u|^{q_0})>N^{q_0}_1\}\cap K_r(\hat{x},\hat{t})|\leq\varepsilon|K_r(\hat{x},\hat{t})|.
\end{eqnarray}
\end{lemma}
\vskip 2mm\noindent{\bf Proof.} From (\ref{d2.20}), there exists a point $(\bar{x},\bar{t})\in Q_{1}\cap K_r(\hat{x},\hat{t})$ such that for any $\rho>0$,
\begin{eqnarray}\label{d2.21}
\fint_{K_{\rho}(\bar{x},\bar{t})\cap Q_7}|D\nabla\times u|^{q_0}dxdt\leq1,
~~\fint_{K_{\rho}(\bar{x},\bar{t})\cap Q_7}|f|^{q_0}dxdt\leq\tilde{\delta}^{q_0},
\end{eqnarray}
\begin{eqnarray}\label{d2.22}
\fint_{K_{\rho}(\bar{x},\bar{t})\cap Q_7}|Dg|^{q_0}dxdt\leq1
~~\mbox{and}~~
\fint_{K_{\rho}(\bar{x},\bar{t})\cap Q_7}|Du-\overline{Du}_{B_{\rho}(\bar{x})}|^{q_0}dxdt\leq\rho^{q_0(1-\gamma)}.
\end{eqnarray}

Let $\widetilde{Q}_{\rho}$ be defined the same as in the proof of Lemma 3.4.
Observe that $\widetilde{Q}_{5r}\subset K_{8r}(\bar{x},\bar{t})\cap Q_7$. Then, (\ref{d2.21}) and (\ref{d2.22}) lead to
\begin{eqnarray}\label{d2.26}
\fint_{\widetilde{Q}_{5r}}|D\nabla\times u|^{q_0}dxdt\leq C,
~~\fint_{\widetilde{Q}_{5r}}|f|^{q_0}dxdt\leq C\tilde{\delta}^{q_0},
\end{eqnarray}
\begin{eqnarray}\label{d2.66}
\fint_{\widetilde{Q}_{5r}}|Dg|^{q_0}dxdt\leq C
~~\mbox{and}~~\fint_{\widetilde{Q}_{5r}}|Du-\overline{Du}_{B_{5r}(\hat{x})}|^{q_0}dxdt\leq Cr^{q_0(1-\gamma)},
\end{eqnarray}
where $C$ depends only on $n$.

Take $\nabla\times$ on both sides of the first equation in (\ref{s0.1}) and then, by $\nabla\times\nabla \pi=0$,
\begin{equation*}
(\nabla\times u)_t-\overline{A}_{B_{4r}(\hat{x})}(t)D^2(\nabla\times u)
=\nabla\times\tilde{f}~~\mbox{in}~~\widetilde{Q}_{4r},\\
\end{equation*}
where $\tilde{f}=\big(A-\overline{A}_{B_{4r}(\hat{x})}(t)\big)D^2u+f$. We decompose $\nabla\times u$ into $w+v$, where $w$ satisfies
\begin{equation*}
\left\{
\begin{array}{ll}
\vspace{2mm} w_t-\overline{A}_{B_{4r}(\hat{x})}(t)D^2w
=\nabla\times\tilde{f}~~\mbox{in}~~\widetilde{Q}_{4r},\\
\vspace{2mm}
w=0~~\mbox{in}~~\partial_p\widetilde{Q}_{4r}.
\end{array}
\right.
\end{equation*}
Then
\begin{equation*}
 v_t-\overline{A}_{B_{4r}(\hat{x})}(t)D^2v
=0~~\mbox{in}~~\widetilde{Q}_{4r}.
\end{equation*}

Instead of (\ref{r.3}), we use (\ref{r.4}) to $w$ to obtain
\begin{equation}\label{d0.63}
\fint_{\widetilde{Q}_{\frac{7r}{2}}}|Dw|^{q_0}dxdt
\leq C\fint_{\widetilde{Q}_{4r}}|\big(A-\overline{A}_{B_{4r}(\hat{x})}(t)\big)D^2u|^{q_0}dxdt
+\fint_{\widetilde{Q}_{4r}}|f|^{q_0}dxdt,
\end{equation}
where $C$ depends only on $n$, $\lambda$ and $p_0$.
(\ref{d2.26}), (\ref{d2.66}) and the classical Calder$\acute{o}$n-Zygmund estimates for (\ref{n.77}) imply
\begin{equation}\label{s7.3}
\begin{aligned}
\fint_{\widetilde{Q}_{4r}}|D^2u|^{q_0}dxdt
\leq&\ C\fint_{\widetilde{Q}_{5r}}\Big(|D\nabla\times u|^{p_0}+\frac{1}{r^{q_0}}|Du-\overline{Du}_{B_{5r}(\hat{x})}|^{q_0}+|Dg|^{q_0}\Big)dxdt\\
\leq&\ C\Big(1+r^{{-q_0}\gamma}\Big),
\end{aligned}
\end{equation}
which is unbounded as $r\rightarrow0$ and where the constant $C$ depends only on $n$ and $p_0$. Then, by (\ref{s7.3}) and $(H_2')$, we obtain
\begin{eqnarray*}
\begin{aligned}
\lefteqn{\fint_{\widetilde{Q}_{4r}}|\big(A-\overline{A}_{B_{4r}(\hat{x})}(t)\big)D^2u|^{q_0}dxdt}\hspace*{24mm}\\
\leq&\ \|A-\overline{A}_{B_{4r}(\hat{x})}(t)\|^{q_0}_{L_\infty(\widetilde{Q}_{4r})}\fint_{\widetilde{Q}_{4r}}|D^2u|^{q_0}dxdt\\
\leq&\ Cr^{{q_0}\gamma}\delta^{q_0}\Big(1+r^{{-q_0}\gamma}\Big)
\leq C\tilde{\delta}^{q_0},
\end{aligned}
\end{eqnarray*}
where $C$ depends only on $n$ and $p_0$. Together with (\ref{d0.63}) and (\ref{d2.26}), it gives
\begin{eqnarray}\label{d2.29}
\fint_{\widetilde{Q}_{\frac{7r}{2}}}|Dw|^{q_0}dxdt\leq C\tilde{\delta}^{q_0},
\end{eqnarray}
where $C$ depends only on $n$, $\lambda$ and $p_0$.

Using Lemma 2.5 (i) to $v$, we get from (\ref{d2.26}) and (\ref{d2.29}) that
\begin{eqnarray}\label{d2.31}
\begin{aligned}
\|Dv\|^{q_0}_{L_\infty(\widetilde{Q}_{3r})}\leq&\ C\fint_{\widetilde{Q}_{\frac{7r}{2}}}|Dv|^{q_0}dxdt
=C\fint_{\widetilde{Q}_{\frac{7r}{2}}}|D\nabla\times u-Dw|^{q_0}dxdt\\
\leq&\ C\Big(\fint_{\widetilde{Q}_{\frac{7r}{2}}}|D\nabla\times u|^{q_0}dxdt+\fint_{\widetilde{Q}_{\frac{7r}{2}}}|Dw|^{q_0}dxdt\Big)
\leq N^{q_0}_0,
\end{aligned}
\end{eqnarray}
where $N_0$ is a constant depending only on $n$, $\lambda$ and $p_0$.

Let $N^{q_0}_1=\max\{2^{q_0}N^{q_0}_0,4^{n+2}\}$.
Then as to derive (\ref{n.63}), we obtain
\begin{equation}\label{d2.6}
\begin{aligned}
\lefteqn{\{(x,t)\in Q_1:\mathcal{M}(|D\nabla\times u|^{q_0})>N^{q_0}_1\}\cap K_r(\hat{x},\hat{t})}\hspace*{16mm}\\
&\subset\{(x,t)\in Q_1:\mathcal{M}_{\widetilde{Q}_{\frac{7r}{2}}}(|Dw|^{q_0})>N^{q_0}_0\}\cap K_r(\hat{x},\hat{t}).
\end{aligned}
\end{equation}
In view of Lemma 2.3 (ii), (\ref{d2.6}) and (\ref{d2.29}), we have
\begin{eqnarray*}
\begin{aligned}
\lefteqn{|\{(x,t)\in Q_1:\mathcal{M}(|D\nabla\times u|^{q_0})>N^{q_0}_1\}\cap K_r(\hat{x},\hat{t})|}\hspace*{12mm}\\
&\leq |\{(x,t)\in Q_1:\mathcal{M}_{\widetilde{Q}_{\frac{7r}{2}}}(|Dw|^{q_0})>N_0\}\cap K_r(\hat{x},\hat{t})|\\
&\mbox{}\hskip1cm\leq \frac{C}{N_0}\int_{\widetilde{Q}_{\frac{7r}{2}}}|Dw|^{q_0}dxdt
\leq \frac{C}{N_0}\tilde{\delta}^{q_0}|\widetilde{Q}_{\frac{7r}{2}}|,
\end{aligned}
\end{eqnarray*}
where $C$ depends only on $n$, $\lambda$ and $p_0$.
Choosing $\tilde{\delta}$ sufficiently small, the above inequality implies
(\ref{d9.20}), which completes our proof.
\qed

\bigskip

\noindent{\bf Proof of Proposition 4.2.}\quad
By Lemma 2.4 and 4.4,
we deduce Proposition 4.2.
\qed

\bigskip

Applying Proposition 4.2 to
$\frac{u}{N^{k-1}_1}$,~~$\frac{\pi}{N^{k-1}_1}$,~~$\frac{f}{N^{k-1}_1}$~~\mbox{and} ~~$\frac{g}{N^{k-1}_1}$,
we have the following corollary.

\begin{corollary}
Under the assumptions of Proposition 4.2,
we have
\begin{eqnarray}\label{d1.41}
\begin{aligned}
\lefteqn{|\{(x,t)\in Q_1:\mathcal{M}(|D\nabla\times u|^{q_0})>N^{q_0k}_1\}|}\hspace*{1mm}\\
\leq&\ 10^{n+2}\varepsilon\Big(|\{(x,t)\in Q_1:\mathcal{M}(|D\nabla\times u|^{q_0})>N^{q_0(k-1)}_1\}|\\
&\mbox{}\hskip0.5cm  +
|\{(x,t)\in Q_1:\Theta>N^{q_0(k-1)}_1\}|+|\{(x,t)\in Q_1:\mathcal{M}(|f|^{q_0})>\tilde{\delta}^{q_0}N^{q_0(k-1)}_1\}|\\
&\mbox{}\hskip1.5cm +|\{(x,t)\in Q_1:\mathcal{M}(|Dg|^{q_0})>N^{q_0(k-1)}_1\}|\Big)
\end{aligned}
\end{eqnarray}
for any $k\in\mathbb{N}_+$.
\end{corollary}

\bigskip

In order to design our iteration to improve the regularity of $D^2u$, we
define
\begin{eqnarray}\label{d1.50}
p_{j+1}=\frac{pp_{j}}{p(1-\frac{\gamma}{2})+p_{j}\frac{\gamma}{2}}~~\mbox{for}~~j=0,1,2,\ldots,
\end{eqnarray}
where $1<p_0<p$ are given in Theorem 4.1 and $\gamma>0$ is the constant in hypothesis ($H_2'$). We can check that $\{p_{j}\}$ is increasing and $\lim\limits_{j\rightarrow\infty}p_j=p$.

\bigskip

We need the following lemma.

\begin{lemma}
Let $1<p_0<p<\infty$, $q_0=\frac{1+p_0}{2}$, and $\Theta$ and $\{p_{j}\}$ be given by (\ref{d2.37}) and (\ref{d1.50}) respectively. Suppose $u\in L_{p}(Q_7;\mathbb{R}^{n})$ and $D^2u\in L_{p_j}(Q_7;\mathbb{R}^{n^3})$. Then
we have $\Theta\in L_{p_{j+1}/q_0}(Q_1)$ and
\begin{eqnarray}\label{d06.48}
\|\Theta\|^{p_{j+1}/q_0}_{L_{p_{j+1}/q_0}(Q_{1})}
\leq C\Big(\|\mathcal{M}(|D^2u|^{q_0})\|^{p_{j}/q_0}_{L_{p_{j}/q_0}(Q_{1})}
+\|\mathcal{M}(|u|^{q_0})\|^{p/q_0}_{L_{p/q_0}(Q_{1})}\Big)
\end{eqnarray}
where $C$ depends only on $n$, $p_0$ and $p$.
\end{lemma}
\vskip 2mm\noindent{\bf Proof.}\quad
From (\ref{d2.37}) and Poincar\'{e}'s inequality, we have
\begin{eqnarray*}
\begin{aligned}
\Theta
=&\  \sup_{r>0}\Bigg\{\Big(r^{-q_0}\fint_{K_r(x,t)\cap Q_7}|Du-\overline{Du}_{B_r(x)}|^{q_0}dyds\Big)^{1-\gamma}\\
&\mbox{}\hskip3cm \Big(\fint_{K_r(x,t)\cap Q_7}|Du-\overline{Du}_{B_r(x)}|^{q_0}dyds\Big)^{\gamma}\Bigg\}\\
\leq&\ C\sup_{r>0}\Bigg\{\Big(\fint_{K_r(x,t)\cap Q_7}|D^2u|^{q_0}dyds\Big)^{1-\gamma}
\Big(\fint_{K_r(x,t)\cap Q_7}|Du|^{q_0}dyds\Big)^{\gamma}\Bigg\}\\
\leq&\ C\big(\mathcal{M}(|D^2u|^{q_0})\big)^{1-\gamma}
\big(\mathcal{M}(|Du|^{q_0})\big)^{\gamma},
\end{aligned}
\end{eqnarray*}
where $C$ depends only on $n$ and $p_0$. Let $\frac{1}{p^*_j}=\frac{1}{2}(\frac{1}{p_j}+\frac{1}{p})$ and then by (\ref{d1.50}), $\frac{p_{j+1}(1-\gamma)}{p_j}+\frac{p_{j+1}\gamma}{p^*_j}=1$. From Young's inequality,
\begin{eqnarray}\label{n.90}
\begin{aligned}
\int_{Q_1}|\Theta|^{\frac{p_{j+1}}{q_0}}dxdt
\leq&\ C\int_{Q_1}\Bigg\{\big(\mathcal{M}(|D^2u|^{q_0})\big)^{\frac{p_{j+1}(1-\gamma)}{q_0}}
\big(\mathcal{M}(|Du|^{q_0})\big)^{\frac{p_{j+1}\gamma}{q_0}}\Bigg\}dxdt\\
\leq&\ C\Bigg\{\int_{Q_1}\big(\mathcal{M}(|D^2u|^{q_0})\big)^{\frac{p_{j}}{q_0}}dxdt
+\int_{Q_1}\big(\mathcal{M}(|Du|^{q_0})\big)^{\frac{p^*_j}{q_0}}dxdt\Bigg\},
\end{aligned}
\end{eqnarray}
where $C$ depends only on $n$, $p_0$ and $p$.
Since, by interpolation inequalities,
\begin{eqnarray*}
\begin{aligned}
\fint_{K_r(x,t)\cap Q_7}|Du|^{q_0}dyds
\leq C\Big(\fint_{K_r(x,t)\cap Q_7}|D^2u|^{q_0}dyds\Big)^{\frac{1}{2}}
\Big(\fint_{K_r(x,t)\cap Q_7}|u|^{q_0}dyds\Big)^{\frac{1}{2}}
\end{aligned}
\end{eqnarray*}
for $r>0$, we obtain
\begin{eqnarray*}
\begin{aligned}
\mathcal{M}(|Du|^{q_0})
\leq C\big(\mathcal{M}(|D^2u|^{q_0})\big)^{\frac{1}{2}}\big(\mathcal{M}(|u|^{q_0})\big)^{\frac{1}{2}},
\end{aligned}
\end{eqnarray*}
where $C$ depends only on $n$ and $p_0$.
Therefore, by Young's inequality, we deduce
\begin{eqnarray*}
\begin{aligned}
\int_{Q_1}\big(\mathcal{M}(|Du|^{q_0})\big)^{\frac{p^*_j}{q_0}}dxdt
\leq&\ C\int_{Q_{1}}\big(\mathcal{M}(|D^2u|^{q_0})\big)^{\frac{p^*_j}{2q_0}}\big(\mathcal{M}(|u|^{q_0})\big)^{\frac{p^*_j}{2q_0}}dxdt\\
\leq&\ C\Big(\|\mathcal{M}(|D^2u|^{q_0})\|^{p_{j}/q_0}_{L_{p_{j}/q_0}(Q_{1})}+\|\mathcal{M}(|u|^{q_0})\|^{p/q_0}_{L_{p/q_0}(Q_{1})}\Big),
\end{aligned}
\end{eqnarray*}
where $C$ depends only on $n$, $p_0$ and $p$.
Substituting it into (\ref{n.90}), (\ref{d06.48}) follows.
\qed

\bigskip

\begin{proposition}
Let $0<\varepsilon<1<p_0<p<\infty$, $q_0=\frac{1+p_0}{2}$, $\{p_{j}\}$ be given by (\ref{d1.50}) and $u$, $\pi$, $f$ and $g$ be given by Theorem 4.1. Assume that $Q_7\subset \Omega_{T}$, $(H_1)$ and $(H_2')$ hold with some $0<\delta<1$ depending only on $\varepsilon$, $n$, $\lambda$, $p_0$ and $p$. If
$D^2u\in L_{p_j}(Q_{7};\mathbb{R}^{n^3})$, then $D^2u\in L_{p_{j+1}}(Q_{\frac{1}{2}};\mathbb{R}^{n^3})$ and
\begin{eqnarray}\label{d6.48}
\begin{aligned}
\|D^2u\|^{p_{j+1}}_{L_{p_{j+1}}(Q_{\frac{1}{2}})}
\leq C\varepsilon\Big(\|D^2u\|^{p_{j}}_{L_{p_{j}}(Q_{7})}
+\|D^2u\|^{p_{j+1}}_{L_{p_{j}}(Q_{7})}\Big)+\tilde{C}M,
\end{aligned}
\end{eqnarray}
where $C$ depends only on $n$, $\lambda$, $p_0$ and $p$, $\tilde{C}$ depends only on $\varepsilon$, $n$, $\lambda$, $p_0$ and $p$ and
\begin{eqnarray*}
\begin{aligned}
M_j=\|u\|^{p}_{L_{p}(Q_{7})}
+\|u\|^{p_{j+1}}_{L_{p}(Q_{7})}+\|f\|^{p_{j+1}}_{L_{p}(Q_{7})}
+\|Dg\|^{p_{j+1}}_{L_{p}(Q_{7})}.
\end{aligned}
\end{eqnarray*}
\end{proposition}
\noindent{\bf Proof.}\quad
Let $N_1$ and $\tilde{\delta}$ be given by Proposition 4.2. Take $\delta=\varepsilon\tilde{\delta}$.
Denote
\begin{eqnarray}\label{n.40}
\tilde{u}=\frac{u}{N}, ~~\tilde{\pi}=\frac{\pi}{N},~~\tilde{f}=\frac{ f}{N}~~\mbox{and}~~\tilde{g}=\frac{g}{N},
\end{eqnarray}
where
\begin{eqnarray*}
N=\bar{C}^{\frac{1}{q_0}}N_1^{-1}\varepsilon^{-\frac{1}{q_0}}\Big(\|u\|_{L_p(Q_{7})}+\|f\|_{L_p(Q_{7})}+\|Dg\|_{L_p(Q_{7})}+\delta\|D^2u\|_{L_{q_0}(Q_{7})}\Big)
\end{eqnarray*}
with the constant $\bar{C}>0$ given by the following (\ref{n1.80}) and depending only on $n$, $\lambda$, $p_0$ and $p$.
It is easy to see that $\tilde{u}$, $\tilde{\pi}$, $\tilde{f}$ and $\tilde{g}$ also satisfy (\ref{s0.1}). From Lemma 2.3 (ii), it follows that
\begin{eqnarray}\label{n.80}
\begin{aligned}
|\{(x,t)\in Q_1:\mathcal{M}(|D\nabla\times \tilde{u}|^{q_0})>N^{q_0}_1\}|
\leq \frac{C}{N^{q_0}_1}\|D\nabla\times \tilde{u}\|^{q_0}_{L_{q_0}(Q_{2})}.
\end{aligned}
\end{eqnarray}
Similar to estimating $\|\nabla\times \tilde{u}\|_{L_{q_0}(Q_{2})}$ in the proof of Proposition 3.7, from (\ref{r.2}) and $(H_2')$, we obtain
\begin{equation*}
\begin{aligned}
\|D\nabla\times \tilde{u}\|_{L_{q_0}(Q_{2})}
\leq C\Big(\delta\|D^2\tilde{u}\|_{L_{q_0}(Q_{7})}++\|\tilde{f}\|_{L_{q_0}(Q_{7})}+\|\tilde{u}\|_{L_{q_0}(Q_{7})}\Big),
\end{aligned}
\end{equation*}
where $C$ depends only on $n$, $\lambda$ and $p_0$.
This and (\ref{n.80}) lead to
\begin{eqnarray}\label{n1.80}
\begin{aligned}
\lefteqn{|\{(x,t)\in Q_1:\mathcal{M}(|D\nabla\times \tilde{u}|^{q_0})>N^{q_0}_1\}|}\hspace*{26mm}\\
\leq&\ \frac{C}{N^{q_0}_1}\Big(\delta\|D^2\tilde{u}\|_{L_{q_0}(Q_{7})}+\|\tilde{f}\|_{L_{q_0}(Q_{7})}+\|\tilde{u}\|_{L_{q_0}(Q_{7})}\Big)^{q_0}\\
\leq&\ \frac{\bar{C}}{(N_1N)^{q_0}}\Big(\delta\|D^2u\|_{L_{q_0}(Q_{7})}+\|f\|_{L_{p}(Q_{7})}
+\|u\|_{L_{p}(Q_{7})}\Big)^{q_0}
\leq \varepsilon|Q_1|.
\end{aligned}
\end{eqnarray}

Applying Corollary 4.5 to $\tilde{u}$, $\tilde{p}$, $\tilde{f}$ and $\tilde{g}$, we have
\begin{equation*}
\begin{aligned}
\lefteqn{\sum^\infty_{k=1}N^{kp_{j+1}}_1|\{(x,t)\in Q_1:\mathcal{M}(|D\nabla\times \tilde{u}|^{q_0})>N^{q_0k}_1\}|}\hspace*{8mm}\\
\leq&\ 10^{n+2}\varepsilon N^{p_{j+1}}_1\sum^\infty_{k=1}\Big(N^{(k-1)p_{j+1}}_1|\{(x,t)\in Q_1:\mathcal{M}(|D\nabla\times \tilde{u}|^{q_0})>N^{q_0(k-1)}_1\}|\\
&\mbox{}\hskip1cm +N^{(k-1)p_{j+1}}_1|\{(x,t)\in Q_1:\tilde{\Theta}>N^{q_0(k-1)}_1\}|\\
&\mbox{}\hskip2cm +N^{(k-1)p_{j+1}}_1|\{(x,t)\in Q_1:\frac{1}{\tilde{\delta}^{q_0}}\mathcal{M}(|\tilde{f}|^{q_0})>N^{q_0(k-1)}_1\}|\\
&\mbox{}\hskip3cm +N^{(k-1)p_{j+1}}_1|\{(x,t)\in Q_1:\mathcal{M}(|D\tilde{g}|^{q_0})>N^{q_0(k-1)}_1\}|\Big),
\end{aligned}
\end{equation*}
where $\tilde{\Theta}$ is defined by (\ref{d2.37}) with $u$ replaced by $\tilde{u}$.
Observe that Lemma 4.6 implies $\tilde{\Theta}\in L_{p_{j+1}/q_0}(Q_{1})$.
Then, by the above inequality, Lemma 2.1, Lemma 2.3 (i) and (\ref{n.40}),
\begin{eqnarray}\label{n.41}
\begin{aligned}
\lefteqn{\|\mathcal{M}(|D\nabla\times \tilde{u}|^{q_0})\|^{p_{j+1}/q_0}_{L_{p_{j+1}/q_0}(Q_{1})}}\hspace*{18mm}\\
\leq&\ C\varepsilon\Big(\|\mathcal{M}(|D\nabla\times \tilde{u}|^{q_0})\|^{p_{j+1}/q_0}_{L_{p_{j+1}/q_0}(Q_{1})}
+\|\tilde{\Theta}\|^{p_{j+1}/q_0}_{L_{p_{j+1}/q_0}(Q_{1})}\Big)+\frac{C\varepsilon^{\frac{p_{j+1}}{q_0}}}{\tilde{\delta}^{p_{j+1}}},
\end{aligned}
\end{eqnarray}
where we use the fact
\begin{eqnarray*}
\begin{aligned}
\|\frac{1}{\tilde{\delta}^{q_0}}\mathcal{M}(|\tilde{f}|^{q_0})\|^{p_{j+1}/q_0}_{L_{p_{j+1}/q_0}(Q_{1})}
\leq \frac{C}{\tilde{\delta}^{p_{j+1}}}\|\tilde{f}\|^{p_{j+1}}_{L_{p_{j+1}}(Q_{7})}
\leq \frac{C\varepsilon^{\frac{p_{j+1}}{q_0}}}{\tilde{\delta}^{p_{j+1}}},
\end{aligned}
\end{eqnarray*}
\begin{eqnarray*}
\begin{aligned}
\|\mathcal{M}(|D\tilde{g}|^{q_0})\|^{p_{j+1}/q_0}_{L_{p_{j+1}/q_0}(Q_{1})}
\leq C\|D\tilde{g}\|^{p_{j+1}}_{L_{p_{j+1}}(Q_{7})}
\leq C\varepsilon^{\frac{p_{j+1}}{q_0}},
\end{aligned}
\end{eqnarray*}
and $C$ depends only on $n$, $\lambda$, $p_0$ and $p$. Take $\varepsilon\leq\frac{1}{2C}$ with the same $C$ in (\ref{n.41}) and then
\begin{eqnarray*}
\begin{aligned}
\|\mathcal{M}(|D\nabla\times \tilde{u}|^{q_0})\|^{p_{j+1}/q_0}_{L_{p_{j+1}/q_0}(Q_{1})}
\leq C\varepsilon\|\tilde{\Theta}\|^{p_{j+1}/q_0}_{L_{p_{j+1}/q_0}(Q_{1})}
+\frac{C\varepsilon^{\frac{p_{j+1}}{q_0}}}{\tilde{\delta}^{p_{j+1}}},
\end{aligned}
\end{eqnarray*}
where $C$ depends only on $n$, $\lambda$, $p_0$ and $p$.
Then
\begin{eqnarray*}
\begin{aligned}
\|D\nabla\times u\|^{p_{j+1}}_{L_{p_{j+1}}(Q_{1})}
\leq&\ C\varepsilon\|\Theta\|^{p_{j+1}/q_0}_{L_{q_{j+1}/q_0}(Q_{1})}+\frac{C}{\tilde{\delta}^{p_{j+1}}}\Big(\|u\|^{p_{j+1}}_{L_{p}(Q_{7})}+\|f\|^{p_{j+1}}_{L_{p}(Q_{7})}\\
&\mbox{}\hskip4.5cm +\|Dg\|^{p_{j+1}}_{L_{p}(Q_{7})}+\delta^{p_{j+1}}\|D^2u\|^{p_{j+1}}_{L_{q_0}(Q_{7})}\Big).
\end{aligned}
\end{eqnarray*}
Recall $\delta=\varepsilon\tilde{\delta}$. Combining this, (\ref{d06.48}) and Lemma 2.3 (i), we get
\begin{eqnarray}\label{n.42}
\begin{aligned}
\|D\nabla\times u\|^{p_{j+1}}_{L_{p_{j+1}}(Q_{1})}
\leq&\ C\varepsilon\Big(\|D^2u\|^{p_{j}}_{L_{p_{j}}(Q_{7})}+\|D^2u\|^{p_{j+1}}_{L_{p_{j}}(Q_{7})}\Big)\\
& +\tilde{C}\Big(\|u\|^{p}_{L_{p}(Q_{7})}+\|u\|^{p_{j+1}}_{L_{p}(Q_{7})}+\|f\|^{p_{j+1}}_{L_{p}(Q_{7})}
+\|Dg\|^{p_{j+1}}_{L_{p}(Q_{7})}\Big),
\end{aligned}
\end{eqnarray}
where $C$ depends only on $n$, $\lambda$, $p_0$ and $p$ and $\tilde{C}$ depends only on $\varepsilon$, $n$, $\lambda$, $p_0$ and $p$.

Using $W^{2,p}$ estimates to elliptic equations (\ref{n.77}), one has
\begin{eqnarray*}
\begin{aligned}
\|D^2u\|^{p_{j+1}}_{L_{p_{j+1}}(Q_{\frac{1}{2}})}
\leq&\ C\Big(\|D\nabla\times u\|^{p_{j+1}}_{L_{p_{j+1}}(Q_{1})}
+\|u\|^{p_{j+1}}_{L_{p_{j+1}}(Q_{1})}
+\|Dg\|^{p_{j+1}}_{L_{p_{j+1}}(Q_{1})}\Big),
\end{aligned}
\end{eqnarray*}
where $C$ depends only on $n$, $p_0$ and $p$.
By substituting (\ref{n.42}) into the above inequality, we obtain
\begin{eqnarray*}
\begin{aligned}
\|D^2u\|^{p_{j+1}}_{L_{p_{j+1}}(Q_{\frac{1}{2}})}
\leq&\ C\varepsilon\Big(\|D^2u\|^{p_{j}}_{L_{p_{j}}(Q_{7})}+\|D^2u\|^{p_{j+1}}_{L_{p_{j}}(Q_{7})}\Big)\\
&+\tilde{C}\Big(\|u\|^{p}_{L_{p}(Q_{7})}
+\|u\|^{p_{j+1}}_{L_{p}(Q_{7})}+\|f\|^{p_{j+1}}_{L_{p}(Q_{7})}
+\|Dg\|^{p_{j+1}}_{L_{p}(Q_{7})}\Big),
\end{aligned}
\end{eqnarray*}
where $C$ depends only on $n$, $\lambda$, $p_0$ and $p$ and $\tilde{C}$ depends only on $\varepsilon$, $n$, $\lambda$, $p_0$ and $p$.
\qed

\begin{corollary}
Under the assumptions of Proposition 4.7,
if $D^2u\in L_{p_j}(Q_{r};\mathbb{R}^{n^3})$, then $D^2u\in L_{p_{j+1}}(Q_{\tau r};\mathbb{R}^{n^3})$ and
\begin{eqnarray}\label{n.145}
\begin{aligned}
\|D^2u\|^{p_{j+1}}_{L_{p_{j+1}}(Q_{\tau r})}
\leq\frac{C\varepsilon}{(r-\tau r)^{(n+2)(\frac{p_{j+1}}{p_j}-1)}}\|D^2u\|^{p_{j+1}}_{L_{p_{j}}(Q_{r})}
+\frac{\tilde{C}}{(r-\tau r)^{2p}}M,
\end{aligned}
\end{eqnarray}
where
$0<\tau<1$, $0< r\leq7$, $C$ depends only on $n$, $\lambda$, $p_0$ and $p$, $\tilde{C}$ depends only on $\varepsilon$, $n$, $\lambda$, $p_0$ and $p$, and
\begin{eqnarray}\label{n.45}
M=\|u\|^{p}_{L_{p}(Q_{7})}+\|f\|^{p}_{L_{p}(Q_{7})}
+\|Dg\|^{p}_{L_{p}(Q_{7})}+1.
\end{eqnarray}
\end{corollary}
\vskip 2mm\noindent{\bf Proof.}\quad
Proposition 4.7 and standard scaling analysis imply
\begin{eqnarray*}
\begin{aligned}
\|D^2u\|^{p_{j+1}}_{L_{p_{j+1}}(Q_{\frac{r}{2}})}
\leq&\ C\varepsilon\Big(\|D^2u\|^{p_{j}}_{L_{p_{j}}(Q_{7r})}
+\frac{1}{r^{(n+2)(\frac{p_{j+1}}{p_j}-1)}}\|D^2u\|^{p_{j+1}}_{L_{p_{j}}(Q_{7r})}\Big)+\frac{\tilde{C}}{r^{2p}}M_j\\
\leq&\
\frac{C\varepsilon}{r^{(n+2)(\frac{p_{j+1}}{p_j}-1)}}\|D^2u\|^{p_{j+1}}_{L_{p_{j}}(Q_{7r})}+\frac{\tilde{C}}{r^{2p}}M
\end{aligned}
\end{eqnarray*}
for $0< r\leq1$, where the last step follows from $\|D^2u\|^{p_{j}}_{L_{p_{j}}(Q_{7r})}\leq\|D^2u\|^{p_{j+1}}_{L_{p_{j}}(Q_{7r})}+1~~\mbox{and}~~M_j+1\leq 4M$. Then, (\ref{n.145}) is obtained by using standard scaling and covering arguments again.
\qed

\bigskip

\noindent{\bf Proof of Theorem 4.1.}\quad
For any integer $j_0\geq1$, define
\begin{eqnarray*}
r_j=1-\big(\frac{1}{2}\big)^{j_0-j+1},~~j=0,1,2,...,j_0.
\end{eqnarray*}
It is easy to see that $r_{j_0}=\frac{1}{2}$ and $\frac{1}{2}\leq r_{j+1}<r_j<1$.
In Corollary 4.8, set $r=r_j$ and $\tau=\frac{r_{j+1}}{r_j}$, we have
\begin{eqnarray*}
r-\tau r=r_j-r_{j+1}=2^{-(j_0-j+1)}.
\end{eqnarray*}
From (\ref{d1.50}),
$\frac{p_{j+1}}{p_j}-1\leq \frac{p}{p_0}\big(1-\frac{\gamma}{2}\big)^j$ and then
\begin{eqnarray*}
(r-\tau r)^{-(n+2)(\frac{p_{j+1}}{p_j}-1)}=(r_j-r_{j+1})^{-(n+2)(\frac{p_{j+1}}{p_j}-1)}\leq 2^{C_0(j_0-j+1)(1-\frac{\gamma}{2})^j}
\end{eqnarray*}
with $C_0=(n+2)\frac{p}{p_0}$.
Then, if $D^2u\in L_{p_j}(Q_{r_j};\mathbb{R}^{n^3})$, by Corollary 4.8, we get $D^2u\in L_{p_{j+1}}(Q_{r_{j+1}};\mathbb{R}^{n^3})$ and
\begin{eqnarray}\label{n6.49}
\begin{aligned}
\|D^2u\|^{p_{j+1}}_{L_{p_{j+1}}(Q_{r_{j+1}})}
\leq C\varepsilon2^{C_0(j_0-j+1)(1-\frac{\gamma}{2})^j}\|D^2u\|^{p_{j+1}}_{L_{p_{j}}(Q_{r_j})}
+2^{2(j_0-j+1)p}\tilde{C}M
\end{aligned}
\end{eqnarray}
for $j=0,1,2,...,j_0-1$, where $C$ depends only on $n$, $\lambda$, $p_0$ and $p$, $\tilde{C}$ depends only on $\varepsilon$, $n$, $\lambda$, $p_0$ and $p$, and $M$ is given by (\ref{n.45}).

Using (\ref{n6.49}) and H\"{o}lder's inequality iteratively from $j=j_0-1$ to $j=0$, it follows that
\begin{eqnarray*}
\begin{aligned}
\|D^2u\|^{p_{j_0}}_{L_{p_{j_0}}(Q_{r_{j_0}})}
\leq&\ 2^{\alpha_{0}+C_0\beta_{0}}
(C\varepsilon)^{\zeta_{0}+1}\|D^2u\|^{p_{j_0}}_{L_{p_{0}}(Q_{r_0})}\\
&\ +
\sum_{j=1}^{j_0}2^{\alpha_{j}+C_0\beta_{j}+\frac{2(j_0-j+2)pp_{j_0}}{p_{j}}}
(C\varepsilon)^{\zeta_{j}}
(\tilde{C}M)^\frac{p_{j_0}}{p_{j}},
\end{aligned}
\end{eqnarray*}
where
\begin{eqnarray*}
\zeta_{0}=\sum_{j=1}^{j_0-1}\frac{p_{j_0}}{p_{j}},
~~\zeta_{j}=\sum_{k=j+1}^{j_0}\frac{p_{j_0}}{p_{k}},
\end{eqnarray*}
\begin{eqnarray*}
\alpha_{0}=\sum_{j=1}^{j_0-1}\Big(\frac{p_{j_0}}{p_{j}}-1\Big),
~~\alpha_{j}=\sum_{k=j}^{j_0-1}\Big(\frac{p_{j_0}}{p_{k}}-1\Big),
\end{eqnarray*}
\begin{eqnarray*}
\beta_{0}=\sum_{j=1}^{j_0}(j_0-j+2)\big(1-\frac{\gamma}{2}\big)^{j-1}\frac{p_{j_0}}{p_{j}}
~~\mbox{and}~~\beta_{j}=\sum_{k=j}^{j_0-1}(j_0-k+1)\big(1-\frac{\gamma}{2}\big)^{k}\frac{p_{j_0}}{p_{k+1}}
\end{eqnarray*}
for $j=1,2,...,j_0$. Note that $D^2u\in L_{p_0}(Q_7;\mathbb{R}^{n^3})$.

Similar to the calculation in the proof of Theorem 3.1, by choosing $\varepsilon$ small enough such that $2^{\frac{4p^2}{p_0}+\frac{4C_0p}{\gamma p_0}}C\varepsilon<1$, we finally obtain
\begin{eqnarray*}
\begin{aligned}
\|D^2u\|^{p_{j_0}}_{L_{p_{j_0}}(Q_{r_{j_0}})}
\leq\big(2^{\frac{p}{p_{0}}+\frac{4C_0p}{\gamma p_0}}C\varepsilon\big)^{j_0}\|D^2u\|^{p_{j_0}}_{L_{p_{0}}(Q_{1})}+
\frac{16^{\frac{p^2}{p_0}}(\tilde{C}M)^\frac{p}{p_{0}}}{1-2^{\frac{4p^2}{p_0}+\frac{4C_0p}{\gamma p_0}}C\varepsilon}.
\end{aligned}
\end{eqnarray*}

Observe $\lim\limits_{j_0\rightarrow\infty}p_{j_0}=p$ and $2^{\frac{p}{p_{0}}+\frac{4C_0p}{\gamma p_0}}C\varepsilon<1$. Letting $j_0\rightarrow\infty$, the above inequality and (\ref{n.45}) lead to
\begin{eqnarray}\label{n.49}
\begin{aligned}
\|D^2u\|_{L_{p}(Q_{\frac{1}{2}})}
\leq&\ C\Big(\big(\|u\|_{L_{p}(Q_{7})}+\|f\|_{L_{p}(Q_{7})}
+\|Dg\|_{L_{p}(Q_{7})}\big)^{\frac{p}{p_0}}+1\Big),
\end{aligned}
\end{eqnarray}
where $C$ depends only on $n$, $\lambda$, $\gamma$, $p_0$ and $p$.

By the same scaling arguments as to derive (\ref{n1}) from (\ref{n.39}), we modify (\ref{n.49}) and obtain
(\ref{n2}).
\qed

\bigskip

\centerline{REFERENCES}
\bibliographystyle{plain}

\end{document}